\documentclass[brochure,12pt]{bourbaki}

\usepackage[applemac]{inputenc}

\usepackage{amssymb}
\usepackage{manfnt}

\usepackage{graphicx}

\usepackage[french]{babel}
\usepackage[T1]{fontenc}
\usepackage{lmodern,textcomp}

\def\n{\noindent}

\def\Q{\mathbf Q}
\def\F{\mathbf F}
\def\P{\mathbf P}
\def\G{\mathbf G}

\def\Z{\mathbf Z}
\def\R{\mathbf R}
\def\C{\mathbf C}
\def\M{\mathbf M}
\def\N{\mathbf N}

\def\deg{\mathop{\rm deg}}
\def\card{\mathop{\rm card}}
\def\dim{\mathop{\rm dim}}
\def\Pol{\mathop{\rm Pol}}
\def\Irr{\mathop{\rm Irr}}
\def\Supp{\mathop{\rm Supp}}

\def\car{\mathop{\rm car}}

\def\Pol{\mathop{\rm Pol}}

\def\Im{\mathop{\rm Im}}
\def\Re{\mathop{\rm Re}}
\def\cp{\mathop{\rm cap}}
\def\Ext{\mathop{\rm Ext}}

\def\l{\ell}
\def\n{\noindent}

\def\zbar{\overline{z}}

\def\A{\mathcal{A}}

\def\sm{\mathchoice
    {\mathbin{\vrule height .72ex width 1.61ex depth -.38ex}}
    {\mathbin{\vrule height .72ex width 1.61ex depth -.38ex}}
    {\mathbin{\vrule height .50ex width 0.85ex depth -.28ex}}
    {\mathbin{\vrule height .20ex width 0.570ex depth -.24ex}}
}
\def\sm{\raisebox{2.20pt}{~\rule{7pt}{1.4pt}~}}
 \def\sm{\mathbin{\raisebox{2.20pt}{\rule{6pt}{1.2pt}}}}

 \date{Mars 2018}
\bbkannee{70\textsuperscript{e} ann\'ee, 2017--2018}
\bbknumero{1146}
\title{Distribution asymptotique des valeurs propres  \\
des endomorphismes de Frobenius}
\subtitle{d'après Abel, Chebyshev, Robinson,...}
\author{Jean-Pierre SERRE}
\address{Collège de France \\
3 rue d'Ulm \\
75005 Paris}
\email{serre@noos.fr}

\begin{document}

\maketitle

 \medskip
\bgroup
\itshape\small
\raggedleft
\`{A} la mémoire de mon vieil ami Michel Raynaud
\par
\egroup

%
\section*{Introduction}

        Lorsque l'on a une suite d'opérateurs linéaires de rang tendant vers l'infini, on peut s'intéresser à la distribution asymptotique de leurs valeurs propres. C'est ce que nous allons faire ici, dans le cas des endomorphismes de Frobenius des variétés abéliennes (autrement dit, des motifs purs de poids 1).
        
        \smallskip 
        
        Fixons un corps fini $\F$, à $q$ éléments. Soit $\A$ l'ensemble des classes d'isogénie de variétés abéliennes sur $\F$ de dimension $>0$. Si $A\in~\A$, soit $P_A(X)$ le polynôme caractéristique de son endomorphisme de Frobenius. C'est un polynôme unitaire de degré $2 \dim A$, à coefficients dans $\Z$,  dont les racines complexes appartiennent au cercle~$C$ de centre $0$ et de rayon $q^{1/2}$; réciproquement, d'après Honda-Tate ([Ta 69]), tout polynôme ayant ces propriétés provient (à une puissance près, cf. lemme 1.8.1) d'une variété abélienne sur $\F$. Comment se répartissent les racines de $P_A$ quand $A$ varie ?
  
     \smallskip 
    Pour donner un sens précis à cette question, il est commode d'utiliser le langage des mesures: si $d = 2 \dim A$,
   écrivons $P_A$ sous la forme  $\prod_{i=1}^d (X-z_i)$ et définissons une mesure $\mu_A$ sur $C$ par
   $\mu_A = \frac{1}{d}\sum \delta_{z_i}$, où $\delta_{z_i}$ désigne la mesure de Dirac en $z_i$; c'est une mesure 
   positive de masse 1. Notons $\M^{\rm ab}$ l'ensemble des mesures sur~$C$ qui sont limites (pour la topologie faible de l'espace des mesures, cf. §1.1) d'une suite de $\mu_A$, avec $A \in \A$. Quelles sont les propriétés des mesures $\mu \in \M^{\rm ab}$, et en particulier quels peuvent être leurs supports ? 
   
   \smallskip
   Nous répondrons partiellement à cette dernière question, en montrant
    que {\it le support d'une mesure de $\M^{\rm ab}$ est, soit fini, soit de capacité $\geqslant q^{1/4}$},  la borne
    $q^{1/4}$ étant optimale. Nous donnerons aussi des exemples où le support de la mesure est {\it un ensemble totalement discontinu
    analogue à l'ensemble triadique de Cantor}. La situation est donc très différente de celle où l'on se limite à des
    jacobiennes de courbes algébriques : dans ce cas, Tsfasman et Vl\u{a}du\c{t} ont montré (cf. [TV 97] et [Se 97])
    que l'on n'obtient essentiellement que des mesures à support égal au cercle $C$, qui ont une densité continue ne s'annulant qu'en un nombre fini de points.

       \medskip

    \n   Le texte est formé de deux §§ suivis de deux Appendices.
    
    \smallskip
    Le §1 contient les démonstrations des énoncés ci-dessus. Il donne d'abord (§1.2) des théorèmes de structure sur
    les éléments de $\M^{\rm ab}$, dans le cadre plus général des entiers algébriques qui sont « totalement » dans un compact fixé de $\C$ (pas nécessairement un cercle). Les démonstrations (§1.4) reposent sur
     la positivité de certaines intégrales du type $\int \log|Q(x)| \mu(x)$, où $Q$ est un polynôme à coefficients dans $\Z$. 
      Ceci fait, il n'y a plus qu'à appliquer un théorème de Robinson ([Ro 64]) pour obtenir les énoncés indiqués plus 
      haut (§1.8).
      
      {\small La plupart des résultats de cette section avaient été obtenus il y a une vingtaine d'années. Je les avais exposés
        à diverses occasions, mais je n'en avais jamais publié les démonstrations. Le présent séminaire me donne l'occasion de combler cette lacune. J'ai appris tout récemment que M.~Tsfasman venait de rédiger un texte ([Ts 18]) qui couvre une partie des §§1.1 à 1.4.\par}

\smallskip      
     Le §2 est consacré à la démonstration du théorème de Robinson, une démonstration très intéressante 
     par les différents arguments qu'elle met en jeu: courbes hyperelliptiques, équation de Pell-Abel et polynômes de Chebyshev.
 
            \smallskip

       Le premier Appendice (Appendice A)
   rassemble quelques définitions et théorèmes standard sur les capacités, tirés principalement des ouvrages de Tsuji
   [Ts 59] et Ransford [Ra 95]. L'Appendice B, rédigé par J. Oesterlé, le complète en faisant le lien avec la théorie
   du potentiel.
   
   \smallskip
   
     Je remercie J-F. Mestre et A. Bogatyrev pour l'aide qu'ils m'ont apportée au sujet du §2. Je remercie aussi J. Oesterlé
     pour sa lecture attentive du manuscrit et ses nombreuses suggestions et corrections.

%
%
%
%

\section{Mesures associées aux entiers algébriques}

%

\subsection{Mesures}
  Les mesures considérées ici sont des mesures de Radon positives sur un espace compact métrisable $K$, au sens de [INT], chap.~III \footnote{Le mode d'exposition de [INT] a été beaucoup critiqué, aussi bien à l'intérieur qu'à l'extérieur de Bourbaki, notamment parce qu'il mêle deux structures différentes : topologie et intégration. Il a cependant l'avantage de bien s'appliquer  aux questions d'équipartition dont il est question ici, car ces questions relèvent à la fois de la topologie et de l'intégration. Le lecteur que cette controverse intéresse pourra consulter l'introduction de [Go~03].}. En d'autres
  termes, ce sont les $\R$-formes linéaires  $f \mapsto \mu(f)$ sur l'espace $C(K)$ des fonctions continues réelles  
  $f$ sur  $K$ telles que:
    
    \smallskip
    (1.1.1) $f \geqslant 0$ \ sur \ $K \  \Rightarrow \ \mu(f) \geqslant 0.$

    \smallskip
    \n On écrit souvent  $\int_K f(x)\mu(x)$ ou $\int_K f\mu$, à la place de
    $\mu(f)$.
    
      \smallskip
    
    La {\it masse} d'une mesure $\mu$ est $\mu(1)$. La {\it mesure de Dirac} en un point $z$ de $K$ est notée~$\delta_z$;
    on a $\delta_z(f) = f(z).$
    
      \smallskip
    Nous aurons besoin plus loin d'intégrer des fonctions semi-continues supérieurement~$F$ sur~$K$ à valeurs dans $\R \cup \{-\infty\}$.
    Par définition (cf. [INT], chap.~IV, §1, appliqué à $-F$ pour transformer 
«supérieurement» en «inférieurement») 
cette intégrale est l'élément de $\R \cup \{-\infty\}$ donné par :
    
      \smallskip
      (1.1.2)  $\int_K F\mu= \inf_{f \geqslant F} \int_Kf\mu$,
      
        \smallskip
     \n  où la borne inférieure porte sur les  $f \in C(K)$ qui majorent $F$.
     
     Par exemple, si  $F$  est la fonction caractéristique d'une partie fermée $T$ de $K$, le nombre $\int_K F\mu$ est la mesure  $\mu(T)$ de $T$ pour $\mu$. Quand $T= \{z\}$ est réduit à un point~$z$, on écrit  $\mu(z)$ à la place
     de $\mu(\{z\})$ ; c'est la {\it masse de $\mu$ en $z$}. Lorsque tous les points sont de masse nulle, on dit que $\mu$ est {\it diffuse} 
(« atomless »), 
cf. [INT], V.5.10.
    
    \medskip
    
    Dans ce qui suit, nous munirons l'espace des mesures de la {\it topologie faible} (également appelée {\it topologie vague}, [INT], chap.~III, §1.9); c'est celle de la convergence simple: une suite $
   (\mu_n)$ de mesures tend vers une mesure $\mu$ si $\mu_n(f) \to \mu(f)$ pour tout $f \in C(K)$.
    Si $F$ est comme ci-dessus,  $\int_K F\mu$ est une fonction semi-continue supérieurement de $\mu$ (puisque c'est une borne inférieure de fonctions continues), autrement dit, on a :
    
      \smallskip
    (1.1.3) $\int_K F\mu \geqslant \limsup \int_K F \mu_n$ \ si \ $\lim \mu_n = \mu$. 
    
      \smallskip
      
      \n En particulier :
      
      (1.1.4) $\int_KF\mu_n \geqslant 0$ pour tout $n \  \Rightarrow \ \int_KF\mu \geqslant 0$.
  
  \smallskip
  L'espace des mesures positives de masse 1 est compact pour la topologie faible, et métrisable, car c'est une partie fermée de la boule unité du dual de l'espace de Banach~$C(K)$, qui est de type dénombrable, cf. [INT], III.1, cor.~3 à la prop.~15.
  
%
\subsection{Entiers algébriques et mesures associées}

À partir de maintenant, $K$ est une partie compacte de $\C$. On va s'intéresser aux entiers algébriques $z \in \C$
dont tous les $\Q$-conjugués appartiennent à $K$, ce que l'on exprime en disant que  $z$  est « {\it totalement dans} $K$ ». 

  De façon plus précise, notons $\Pol_K$ l'ensemble des polynômes unitaires de degré $>0$, à coefficients dans $\Z$, dont toutes les racines appartiennent à $K$. Si $P$
 est un tel polynôme, de degré $d$ et de racines $z_1,\ldots,z_d$, on lui associe la mesure $\delta_P$ sur $K$ définie par:
 
   \medskip

    (1.2.1)  $\delta_P = \frac{1}{d} (\delta_{z_1}+ \cdots + \delta_{z_d})$, \ i.e. \ $\delta_P(f) =\frac{1}{d} \sum f(z_i)$\  pour  tout $f \in C(K).$
    
     \medskip

Notons  $\M$ (ou $\M_K$ lorsque l'on veut préciser $K$) l'adhérence pour la topologie faible de la famille des mesures
$\delta_P$, où $P$ parcourt $\Pol_K$.

 \smallskip
  C'est la structure de $\M$ qui nous intéresse. On a tout d'abord :
  
   \smallskip
   \n {\bf Proposition 1.2.2.}  {\it L'espace $\M$ est convexe et compact.}
   
    \smallskip
   La compacité résulte de ce que  $\M$ est une partie fermée de l'espace des mesures de masse 1 sur  $K$, qui est compact. Pour la convexité, on observe que, si  $P_1, \ldots, P_m$ appartiennent à  $\Pol_K$, il en est de même
   de leurs produits  $P_1^{a_1}\cdots P_m^{a_m}$ avec $a_i \in \N$, et les~$\delta_P$ correspondants sont denses dans le simplexe
   de sommets les $\delta_{P_i}$. 
   
   \smallskip
   
   \n {\bf Corollaire 1.2.3.}  {\it L'espace $\M$ est l'enveloppe convexe fermée de l'ensemble des~$\delta_P$}.
   
   \medskip

   Notons $\Irr_K$ le sous-ensemble de $\Pol_K$ formé des polynômes irréductibles, et soit $I_K$ l'ensemble des $\delta_P, \ P\in \Irr_K$. Les éléments de~$I_K$ sont linéairement indépendants, et leur enveloppe convexe fermée est $\M$. Lorsque $\Irr_K$ est fini, cela donne la structure de~$\M$:
   c'est le simplexe dont l'ensemble des sommets est $I_K$.
   
   \smallskip
   
    Supposons $\Irr_K$ infini. C'est un ensemble dénombrable. Numérotons ses éléments : $P_1, P_2, \ldots$ Pour tout $n \geqslant 1$, soit $\M_n$ l'enveloppe convexe fermée des $\delta_{P_i}, i \geqslant n$. On a :
    
    \smallskip
    (1.2.4) \ $\M = \M_1  \supset \M_2   \supset \cdots$
   
   \smallskip
    Posons :
    
    \smallskip
    
    (1.2.5) $\M_\infty = \bigcap_{n \geqslant 1} \M_n$. 
    
    \smallskip
    C'est un convexe compact non vide; il ne dépend pas de la numérotation des éléments de $\Irr_K$.
   
   \medskip
   
   \n {\bf Théorème 1.2.6}. {\it Soit $\mu \in \M_\infty$, et soit $S = \Supp \mu$ son support. Alors} :

   (1.2.7)  {\it $\mu$ est diffuse} (cf. §1.1). 
      
   (1.2.8) {\it La capacité $\cp(S)$ de $S$ est $\geqslant 1;$  si  elle est égale à $1$, alors $\mu$ est la mesure d'équilibre de $S$.} 
   
   (1.2.9) {\it Si $K \subset \R$, l'ensemble $S$ est réduit au sens de A.4.7.}
   
  (Pour tout ce qui concerne les capacités et les mesures d'équilibre, voir les deux Appendices; les références
correspondantes commencent par les lettres A et B.)  

\medskip
\n {\bf Corollaire 1.2.10} (Fekete [Fe 23], Satz XI). {\it Si $\cp(K) < 1$, alors $\Irr_K$ est fini.}
   
   \medskip
  \n La démonstration du th.~1.2.6 sera donnée au §1.4.4.
   
   \medskip
   La structure de $\M$ se ramène à celle de $\M_\infty$ par le théorème suivant :

   \smallskip
   \n {\bf Théorème 1.2.11.} {\it Soit $\mu \in \M$. Il existe une suite et une seule de nombres réels positifs $c_0, c_1, c_2, \ldots$
   tels que $\sum_{i\geqslant 0} c_i = 1$, et que $:$
 
    \smallskip  
   $(1.2.12)$ \ \  $\mu = \sum_{i \geqslant 1} c_i \delta_{P_i} + \nu$, \  avec} \ $\nu \in c_0\M_\infty$.

     \n {\small Rappelons que $P_1, P_2, \ldots$ sont les différents éléments de $\Irr_K$, numérotés dans un ordre arbitraire.\par}
      
      \smallskip
 \n  La démonstration sera donnée au §1.4.5.

     \smallskip
   \n {\bf Corollaire 1.2.13}. {\it Si $\mu \in \M$ n'est pas combinaison linéaire d'un nombre fini de  $\delta_{P_i}$, la capacité de 
   $\Supp \mu$ est $\geqslant 1$.}

  \smallskip
  
    Si $\nu \neq 0$, cela résulte de (1.2.8) appliqué à $c_0^{-1}\nu$. Si $\nu = 0$,  il y a une infinité de $c_i$ qui sont non nuls, et $\Supp \mu$
    contient les zéros des $P_i$ correspondants; d'après (1.2.10) appliqué à $\Supp \mu$, on a $\cp(\Supp \mu) \ \geqslant 1.$

      \medskip
   \n {\it Remarque}. La somme infinie

      \medskip
   (1.2.14) $\mu_{\rm at} =  \sum_{i \geqslant 1} c_i \delta_{P_i}$

      \medskip
 \n   est une mesure {\it atomique} ([INT], III.1.3). La formule $\mu = \mu_{\rm at} + \nu$ donne la décomposition canonique de $\mu$ en partie atomique et partie diffuse, cf. [INT] V.5.10, prop.~15. 
 
 \medskip

 Les  th.~1.2.6 et 1.2.11 entraînent :
   
      \smallskip
   \n {\bf Théorème 1.2.15.}  {\it Soit $\mu \in \M$.}
   
    (i) {\it \  $\mu \in \M_\infty \Longleftrightarrow \mu$ est diffuse.}

   (ii) {\it La masse $\mu(z)$ de $\mu$ en un point $z\in K$ est $0$ si  $z$ n'est pas un entier algébrique totalement dans $K$.}
   
   (iii) {\it Si $z, z' \in K$ sont deux entiers algébriques conjugués, on a $\mu(z)  = \mu(z')$.}

   \smallskip
   
   \n {\it Remarque sur} (ii). Cette propriété exprime une sorte d'{\it indépendance} des mesures $\delta_{P_i}$, plus forte que la simple indépendance linéaire (le lemme 1.3.1 ci-dessous va dans la même direction). En général, l'enveloppe convexe fermée d'une suite de mesures discrètes contient bien d'autres mesures discrètes. Par exemple, sur $K = [0,1]$, l'enveloppe convexe fermée des mesures $\delta_z$, avec $z \in \Q \cap K$,  est l'ensemble de {\it toutes les mesures positives de masse $1$ sur $K$}, y compris les mesures de Dirac en des points irrationnels.

%
%

\subsection{Lemmes de positivité}

  Les démonstrations des théorèmes du §1.2 sont basées sur la positivité de certaines intégrales portant sur les mesures $\mu \in \M$, les fonctions intégrées étant du type  $z \mapsto \log |Q(z)|$, où  $Q$ est un polynôme  \footnote{L'utilisation d'intégrales  $\int \log|Q(z)| \mu(z)$ est une technique standard depuis l'application que C.~Smyth en a faite pour l'estimation des traces des entiers algébriques totalement positifs, cf. [Sm 84].} à coefficients dans $\Z$. Noter que, lorsque  $Q(z) = 0$, on a $\log |Q(z)| = - \infty$; ainsi $\log |Q|$ est une fonction continue sur $K$ à valeurs
  dans $\R \cup \{- \infty\}$, et on peut lui appliquer le §1.1.2, ce qui donne un sens aux intégrales de type $\mu(\log |Q|)$.

   \smallskip

\n {\bf Lemme 1.3.1.} {\it Soit $P \in \Pol_K$ et soit $Q \in \Z[X]$, $Q \neq 0$. 

On a
$\delta_P(\log |Q|) = - \infty$ si $P$ et~$Q$ ont une racine commune, et $\delta_P(\log |Q|) \geqslant0$ sinon.}
 
    \smallskip
    
 \n    {\it Démonstration}. Soit $d$ le degré de $P$ et soient $z_1,\dots,z_d$ ses racines. On a
    
     \medskip
(1.3.2) \ \ 
$\delta_P(\log |Q|) = \frac{1}{d} \sum \log |Q(z_i)| = \frac{1}{d} \log |\prod Q(z_i)| = \frac{1}{d} \log |R|$,

 \medskip 
 
 \n où $R$ est le résultant de  $P$ et $Q$ ([A IV], §6.6). Comme $R$ est un entier, on a $\log |R| \geqslant 0$
si $R \neq 0$ et $\log |R| = - \infty$ si $R=0$. D'où (1.3.1), puisque $R=0$ si et seulement si $P$ et~$Q$ ont une racine commune.

 \smallskip
 
\n {\bf Lemme 1.3.3.} {\it Soit $\mu \in \M_n, \  n\geqslant 1$, et soit \ $Q \in \Z[X], \ Q \neq 0$. 
Supposons que $Q$ ne soit divisible par aucun $P_i, \ i \geqslant n$. On a alors
 $\mu(\log |Q|) \geqslant 0$.}
 
    \smallskip
 \n    {\it Démonstration}. Lorsque $\mu$ est combinaison linéaire finie des $\delta_{P_i}, i \geqslant~n$, la positivité de
$\mu(\log |Q|) $ résulte de (1.3.1). Le cas général s'en déduit par passage à la limite en utilisant (1.1.4).

 \smallskip

   \n {\bf Lemme 1.3.4.} {\it On a $\mu(\log |Q|) \geqslant 0$ si $\mu \in \M_\infty$ et} $Q \in \Z[X], Q \neq 0$. 
   
\smallskip

  Cela résulte du lemme précédent, appliqué en prenant $n$ assez grand. 
  
\medskip

  Nous allons maintenant nous occuper d'intégrales doubles (il y a des énoncés analogues pour les intégrales triples, etc. --- nous n'en aurons pas besoin).
  \smallskip
  
   \n {\bf Lemme 1.3.5.} {\it Soit $\mu$ une mesure positive sur $K$, soit $z\in K$, et soit $F(x,y)$ une fonction semi-continue supérieurement sur $K \times K$, à valeurs dans $\R \cup \{-\infty\}$. On a}:
   
   \smallskip
   (1.3.6) \ $\iint_{K\times K} F(x,y) \ \delta_z(x)\mu(y) \ = \ \int_K F(z,y) \ \mu(y)$.
  
  \smallskip
  
\n {\it Démonstration.} Lorsque $F$ est continue à valeurs dans $\R$, la formule (1.3.6) est une forme élémentaire du théorème de Lebesgue-Fubini : on la démontre en se ramenant au cas où $F(x,y)$ est de la forme $f(x)g(y)$. Le cas général s'en déduit en écrivant $F$ comme borne inférieure de fonctions continues à valeurs dans $\R$, et en appliquant
 (1.1.2) aux deux membres de (1.3.6).

\n{\it Démonstration.} Cela résulte du théorème de Lebesgue-Fubini ([INT], chap.V, §8, prop.5), appliqué à $ y \mapsto - F(z,y) + C$,
 où  \mbox{$C \geqslant \sup_{y\in X} F(z,y)$} --- le signe « moins » et la constante $C$ sont dus au fait que Bourbaki traite le cas des fonctions semi-continues inférieurement, à valeurs positives.
 
  \medskip

  \n {\bf Lemme 1.3.7.} {\it Soient $\mu,\nu \in \M_\infty$, et soit $Q \in \Z[X,Y], \ Q \neq 0$. On a}:

(1.3.8) \quad $\iint_{K\times K} \log|Q(x,y)| \ \mu(x)\nu(y) \geqslant 0.$ 

\medskip

  \n {\it Démonstration}. Notons $I(\mu,\nu)$ le membre de gauche de (1.3.8). Il a un sens pour tout couple de mesures $\mu,\nu$ sur $K$. 
  
  Choisissons le couple $(\delta_P, \nu)$, avec 
   $P\in \Irr_K$; soit $d$ le degré de~$P$, et soient $z_1,\ldots,z_d$ ses racines.  Supposons que $P(X)$ ne divise pas $Q(X,Y)$, i.e. que les polynômes $Q(z_i,Y)$ soient $\neq 0$, et posons
 $H_P(Y) = \prod Q(z_i,Y)$; c'est un polynôme en $Y$, non nul, et à coefficients entiers. On a :
  
  \medskip
  
  $I(\delta_P,\nu) = \iint_{K\times K} \log|Q(x,y)|\ \delta_P(x)\nu(y)$.
  
  \smallskip
  
   \hspace{14mm} $ = \frac{1}{d} \sum \iint_{K\times K} \log|Q(x,y)| \ \delta_{z_i}(x)\nu(y)$,

  \smallskip
  
    \hspace{14mm} $ = \frac{1}{d} \sum \int_K \log|Q(z_i,y)| \ \nu(y)$,  d'après (1.3.6),
   
   \smallskip
    \hspace{14mm} $ = \frac{1}{d} \int_K \log|H_P(y)| \ \nu(y) \ \geqslant 0$,  d'après le lemme 1.3.4.
   
   \medskip
   On déduit de là que $I(\delta_P,\nu) \geqslant 0$ pour tous les $P\in \Irr_K$, sauf ceux qui divisent $Q(X,Y)$. Par combinaisons linéaires et passages à la limite, on en déduit le même résultat pour
   $I(\mu, \nu)$, avec $\mu \in \M_n$ pour $n$ assez grand. D'où (1.3.8), sous une forme un peu plus précise :
   
   \smallskip
(1.3.9)  \   $I(\mu, \nu) \ \geqslant 0$ {\it si $\nu \in \M_\infty$ et $\mu \in \M_n$ et si aucun des $P_i(X), i \geqslant n$, ne divise}  $Q(X,Y)$. 
    
\medskip
Le cas particulier le plus utile de (1.3.8) est celui où $Q = X-Y$:

\medskip
\n {\bf Corollaire 1.3.10.} {\it Si $\mu,\nu \in \M_\infty$, on a $\iint_{K\times K} \log|x-y| \mu(x)\nu(y) \geqslant~0$.}

%
\subsection{Démonstration du th.~1.2.6 et du th.~1.2.10}

L'énoncé suivant est un cas élémentaire du théorème de décomposition des mesures
en partie atomique et partie diffuse ([INT], V.5.10, prop.~15) :

\smallskip
\n {\bf Proposition 1.4.1.} {\it Soit $\mu$ une mesure positive sur $K$, et soit $c$ sa masse en un point $z\in K$.
Alors la mesure $ \mu - c\delta_z$ est positive.}

\smallskip

\n {\it Démonstration}. Soit $f$ une fonction continue positive sur $K$.
Il nous faut prouver que  $\mu(f) \geqslant cf(z).$ Si $f(z)=1$, cela résulte de la définition de $c$
rappelée au §1.1. Le cas général en résulte par homogénéité.

\smallskip
\n {\bf Corollaire 1.4.2.} {\it Soient $\mu, z, c$ comme ci-dessus, avec $c > 0$. Soit $F$ une fonction semi-continue
supérieurement sur $K$, à valeurs dans $\R \cup \{-\infty\}$, et  telle que  $F(z) = -\infty$. Alors $\mu(F) = -\infty$.}

 \smallskip

\n {\it Démonstration.} D'après la prop.~1.4.1, on a $\mu = c\delta_z + \nu$, où  $\nu$ est une mesure positive; ainsi,
$\nu(F)$ a un sens.
On a $\mu(F) = cF(z) + \nu(F) = -\infty + \nu(F) = -\infty$.

\n  [Rappelons que $-\infty + x = -\infty$ pour tout $x \in \R \cup \{-\infty\}$.]

\smallskip

\n {\bf Corollaire 1.4.3.} {\it Si $\mu \in \M_n, n \geqslant 1$, on a $\mu(z)=0$ pour toute racine $z$ de l'un des polynômes
$P_i, \  i < n.$}

\smallskip

\n {\it Démonstration.} Si l'on avait $\mu(z)>0$ avec $P_i(z) = 0$ et $i<n$, on aurait  $\mu(\log |P_i|) = -\infty$ d'après le cor.~1.4.2, ce qui contredirait le lemme 1.3.3, appliqué à $Q = P_i$.

\medskip

On va maintenant s'occuper des démonstrations des énoncés du §1.2.

\smallskip
\n 1.4.4. {\it Démonstration du th.~1.2.6.} Soit $\mu \in \M_\infty$. Si $\mu$ avait une masse non nulle en un point $z \in K$, la mesure $\mu \otimes \mu$ aurait une masse non nulle au point $(z,z)$ de $K \times K$; comme la fonction $\log|x-y|$
vaut $-\infty$ en $(z,z)$, le cor.~1.4.2 montrerait que l'intégrale $I(\mu)=\iint_{K\times K} \log|x-y| \mu(x)\mu(y)$ est égale à  $-\infty$, ce qui contredirait
le cor.~1.3.10, qui dit que $ I(\mu)\geqslant 0$.  La mesure $\mu$ est donc diffuse. 

Si $S$ est son support, la positivité de $I(\mu)$ et la caractérisation (A.3.2) de $\cp S$ montrent que $\cp(S) \geqslant 1$; le fait que $S$ soit réduit résulte du cor.~A.4.9. Si $\cp(S) = 1$, alors $I(\mu) = 0$ et $\mu$ est la mesure d'équilibre de $S$, vu l'unicité de celle-ci.

\smallskip\smallbreak
\n 1.4.5. {\it Démonstration du th.~1.2.11.}

\smallskip

Commençons par un résultat partiel :

\smallskip\nobreak
\n {\bf Proposition 1.4.6.} {\it Soit  $n > 0$ et soit $\mu \in \M.$}

\smallskip\nobreak
(i) {\it Il existe une famille et une seule de nombres réels positifs  $c_0, c_1,\ldots,c_{n-1}$, avec $\sum_{i\geqslant 0} c_i = 1$ telle que} 
$\mu - \sum_{i>0} c_i\delta_{P_i} \ \in \  c_0\M_n.$ 

\smallskip
(ii) {\it Soit $i$ tel que $0<i<n$, et soit  $d_i = \deg P_i$. On a $c_i = d_i\mu(z)$ pour toute racine $z$ de $P_i$.}

\smallskip

\n {\it Démonstration.} Le sous-espace de $\M$ formé des mesures de la forme :

\smallskip

   $c_1\delta_{P_1} + \cdots + c_{n-1}\delta_{P_{n-1}} + c_0\nu$, \ avec \ $c_i \geqslant 0, \ \sum_{0\leqslant i<n} c_i=1, \ \nu \in \M_n,$
   
   \smallskip
   
  \n est une partie convexe de $\M$ qui est fermée (car image d'un compact) et qui contient tous les $\delta_P, P \in \Irr_K$.
   C'est donc $\M$, ce qui démontre la partie « existence » de (i). 
   
   Si $\mu$ est de la forme $c_1\delta_{P_1} + \cdots + c_{n-1}\delta_{P_{n-1}} + c_0\nu$ comme ci-dessus, et si $z$ est l'une des racines de l'un des $P_i$,  on a
   $\delta_{P_j}(z)=0$ pour $ j \neq i$,\ $ \delta_{P_i} (z)= 1/d_i$ et $\nu(z)=0$, cf. cor.~1.4.3. D'où (ii); l'assertion d'unicité
   de (i) en résulte.
      
  \smallskip
  
  \smallskip
\n {\bf Corollaire 1.4.7.} {\it Si $\mu \in \M$, on a $\mu \in \M_n$ si et seulement si $\mu(z) = 0$ pour toute racine $z$ de  $P_1 \cdots P_{n-1}$.}

\medskip
\n {\it Fin de la démonstration du théorème 1.2.11}.

\smallskip
 Soit $\mu \in \M$. Pour tout $i >0$, choisissons une racine $z$ de $P_i$ et posons $c_i = d_i\mu(z)$, où $d_i = \deg P_i$.  Si $n > 1$, posons  $\gamma_n = 1 - \sum_{i<n} c_i$. D'après la prop.~1.4.6, $\gamma_n$ est $\geqslant 0$, et l'on a :
 
 \smallskip
 (1.4.8)   \   $ \mu = \sum_{i < n} c_i\delta_{P_i} \ + \nu_n$, \ \ avec \ \ $\nu_n \in  \gamma_n \M_n.$
 
   \smallskip
   
 \n  La série 
   
   \smallskip
   
   (1.4.9)  \ $\mu_{\rm at} = \sum c_i\delta_{P_i} $
   
   \smallskip
   
\n    est convergente, non seulement pour la topologie faible, mais aussi pour la topologie forte (celle donnée par la norme de l'espace des mesures). Posons  $\nu = \mu - \mu_{\rm at}$. En comparant (1.4.8) et (1.4.9), on obtient :
   
    \smallskip
   (1.4.10) \ $\nu = \lim_{n \to \infty} \nu_n$,
   
    \smallskip
   \n la limite étant prise au sens de la topologie forte.
   C'est une mesure positive de masse $c_0 = 1 - \sum_i c_i = \lim \gamma_n$. Nous allons voir qu'elle appartient à  $c_0\M_\infty$, ce qui démontrera la partie « existence » du th.~1.2.10 (la partie « unicité » résulte de la prop.~1.4.6).
   C'est clair si $c_0 = 0$, car cela implique que la masse de $\nu$ est 0, donc que $\nu=0$.
   Supposons $c_0 > 0.$ On a 
   
   \smallskip
   
   (1.4.11) $c_0^{-1}\nu = \lim_{n \to \infty} c_0^{-1} \nu_n.$
   
     \smallskip

   Les masses des mesures positives  $(c_0^{-1}-\gamma_n^{-1})\nu_n$ tendent vers $0$ pour $n \to \infty$. On peut donc récrire (1.4.11) sous la forme :
   
     \smallskip

   (1.4.12) $c_0^{-1}\nu = \lim_{n \to \infty} \gamma_n^{-1} \nu_n.$
    
    \smallskip
   On a $\gamma_n^{-1} \nu_n \in \M_n$. Pour tout entier $m$, on a donc $\gamma_n^{-1} \nu_n \in \M_m$ pour
   $n$ assez grand, et (1.4.12) entraîne $c_0^{-1}\nu  \in 
   \M_m$. Comme ceci est vrai pour tout $m$, on a
   $c_0^{-1}\nu \in \M_\infty$, i.e. $\nu \in c_0\M_\infty$, ce qui termine la démonstration du th.~1.2.10.
   
%
%
\subsection{Deux exemples: le cercle unité et l'intervalle $[-2,2]$}
   
   Prenons pour $K$ le cercle unité $C_1$, i.e. l'ensemble des $z\in \C$ de module 1. Les entiers algébriques totalement dans $K$ sont les racines de l'unité, cf. [Kr~57]. Les éléments de~$\Irr_K$ sont les polynômes cyclotomiques $\Phi_n, \ n = 1, 2, \ldots $ La capacité de $K$ est 1; d'après (1.2.8), cela montre que le seul élément de $\M_\infty$ est la mesure d'équilibre~$\mu_K$ de~$K$, qui n'est autre que la mesure de masse~1 invariante par rotation, autrement dit la mesure $\frac{1}{2\pi}d\varphi$, si l'on écrit les éléments de~$K$ comme $z=e^{i\varphi}$. D'après le th.~1.2.10, tout élément~$\mu$ de~$\M$ s'écrit de façon unique sous la forme
   
   \smallskip
   
(1.5.1)   \ $\mu = c_0\mu_K  +  \sum_{n=1}^\infty c_n \delta_{\Phi_n},$

\smallskip

\n où les $c_i$ sont des nombres réels positifs tels que $\sum_{n=0}^\infty c_n= 1$.
   
   On peut montrer que l'application  $\mu  \mapsto  (c_n)_{n\geqslant 1}$ est un homéomorphisme de $\M$ sur le sous-espace du cube infini $[0,1] \times [0,1] \times \cdots
  $ formé des $(c_n)$ tels que $\sum_{n\geqslant 1} c_n  \leqslant 1$.
  
  \smallskip
 {\small
\advance\leftskip 2 \parindent 
\parindent0pt
\leavevmode\hskip-\leftskip\hbox to\leftskip{\includegraphics[scale=.5]{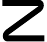}\hfil}%
Noter que l'application $\M \to [0,1]$ donnée
par $\mu \mapsto c_n$ est continue si $n > 0$, mais ne l'est pas si $n=0$ ; la coordonnée $c_0$ ne joue pas le même rôle que les autres. D'ailleurs le sous-espace du cube infini formé des $(c_n)_{n \geqslant 0}$ de somme 1 n'est pas fermé.\par}
   
 \smallskip
 
 Du cercle $C_1$ on passe à l'intervalle $I = [-2,2]$ par $z \mapsto z + \zbar$, cf. A.6. Les entiers algébriques totalement dans $I$ sont les $z + \zbar$, où $z$ est une racine de l'unité. On a une bijection naturelle des
 polynômes irréductibles correspondants, d'où aussi un isomorphisme de~$M_{C_1}$ sur~$M_I$. L'unique élément de
 $\M_{\infty,I}$ est la mesure d'équilibre $\mu_I = \frac{1}{\pi}\frac{dx}{\sqrt{4-x^2}}$. 
   
   \smallskip
   
   \n {\it Remarque.}  Au lieu de prendre  $K=I$, comme nous venons de le faire, on pourrait choisir pour $K$
   n'importe quel ensemble compact contenant $I$ et s'intéresser à l'espace $\M_{K,I} $ formé des {\it $\mu \in \M_K$ dont le support est contenu dans $I$}. A priori, cet espace pourrait être strictement plus grand que $\M_I$ -- cela se produit dans d'autres cas. Mais cela ne se produit pas  ici: on a $\M_{K,I} = \M_I$;  cela se voit en remarquant
   que, comme $\cp(I) = 1$, la seule mesure diffuse de      $\M_{K,I}$ est $\mu_I$. Même chose pour le cercle $C_1$.

%
%
\subsection{Le théorème de Robinson et ses applications} 

  
    
  \smallskip   
    
    \n {\bf Théorème 1.6.1} (Robinson [Ro 64]). {\it Soit $E$ la réunion d'un ensemble fini d'intervalles fermés de $\R$. Supposons $\cp(E) > 1$.
       Il existe une infinité d'entiers algébriques totalement dans $E$.}
     
     \smallskip
       De façon plus précise :
     
     \smallskip
     \n {\bf Théorème 1.6.2}. {\it Il existe une suite de polynômes $P_n \in \Pol_E$ telle que \ $\lim_{n \to \infty} \mu_{_{P_n}} =~\mu_E.$}

     \smallskip
   
     La démonstration des théorèmes 1.6.1 et 1.6.2 fera l'objet du §2.
     
     \medskip
     
     \n {\it Remarque}. 
     
    Le th.1.6.1 a été démontré par Robinson deux ans après qu'il ait traité le cas particulier (déjà très intéressant) où $E$ est un intervalle, cf. [Ro 62]. C'est l'analogue sur  
      $\R$ du théorème suivant, dû à Fekete-Szeg\H o ([FS 55]) :
     
         \smallskip
       {\it Soit $E$ un compact de $\C$ de capacité $\geqslant 1$, invariant par conjugaison. Pour tout voisinage 
       $U$ de $E$ dans $\C$, il existe une infinité d'entiers algébriques qui sont totalement dans $U$.}

       {\small(Noter que, si $E \subset \R$, les entiers algébriques en question ne sont pas nécessairement dans $\R$ ; ils en sont  « arbitrairement voisins », ce qui n'est pas suffisant pour ce qui nous intéresse.) \par}
       
               \medskip   
        
             Soit maintenant $K$ un compact de $\R$, de capacité $> 1$, qui soit réunion d'un nombre fini d'intervalles.
     Nous allons donner deux applications du th.1.6.2 à $\M_\infty= \M_{K,\infty}$. Tout d'abord :

 \smallskip   
    
    \n {\bf Théorème 1.6.3.} {\it Soit $E$ une partie compacte de $K$ de capacité $\geqslant1,$ et soit $\mu_E$ sa mesure d'équilibre. On a $\mu_E \in \M_\infty$.}

\smallskip
\n {\it Démonstration.} Pour tout $\varepsilon \geqslant 0$, soit $E_\varepsilon$ l'ensemble des points de $K$ dont la distance à $E$ est $\leqslant \varepsilon$. On a  $E_0 = E$ et $E_\varepsilon = K$ lorsque $\varepsilon$ est assez grand. De plus, si $\varepsilon > 0$, $E_\varepsilon$ est une réunion finie d'intervalles. Si $\cp(E) > 1$, il en est de même des $E_\varepsilon$ et l'on a $\mu_{E_\varepsilon} \in \M$ d'après le th.1.6.2. Comme $E$
est l'intersection des $E_\varepsilon \ (\varepsilon > 0)$, on a $\mu_E = \lim_{\varepsilon \to 0}\mu_{E_\varepsilon}$ d'après (A.5.2), d'où $\mu_E \in \M$; comme $\mu_E$ est une mesure d'équilibre, elle est diffuse, donc appartient à $\M_\infty$, cf. th.1.2.15. Supposons maintenant $\cp(E)=1$. Soit $\alpha$ la borne supérieure des $\varepsilon \geqslant 0$ tels que $\cp(E_\varepsilon) = 1$. Si $\beta > \alpha$, le th.1.6.2 montre que $\mu_{E_\beta}$ appartient à $\M$, et d'après (A.5.2), il en est de même de $\mu_{E_\alpha}$. Si $\alpha = 0$, on a $E_\alpha = E$, d'où $\mu_E \in \M$. Supposons $\alpha > 0$; alors $E_\alpha$ contient
$E'_\alpha = \cup_{\gamma < \alpha} E_\gamma$. L'ensemble $E_\alpha \sm E'_\alpha$ est l'ensemble des $x\in K$ dont la distance à $E$ est égale à $\alpha$. C'est donc un ensemble fini (plus généralement, si $\Omega$ est un compact de $\R$,
l'ensemble des points de $\R$ dont la distance à $\Omega$ est $
\alpha$ est un ensemble fini); sa capacité est donc~$0$. En appliquant (A.5.5) à $E_\alpha$ on en conclut que $\cp(E_\alpha) = \sup \cp(E_\gamma)=1.$ Les compacts $E$ et $E_\alpha$ ont même capacité, et l'on a $E \subset E_\alpha$. Ils ont donc même mesure d'équilibre:
cela résulte de l'unicité de la mesure d'équilibre de $E_\alpha$.
Comme on a vu que $\mu_{E_\alpha}$ appartient à $\M$, on a $\mu_E \in \M$, d'où $\mu_E \in \M_\infty$. 

   \medskip
   
 \n {\it Exemple.} Le th.1.6.3 s'applique notamment aux intervalles de longueur $\geqslant 4$ contenus dans $K$ : les mesures d'équilibre
 correspondantes appartiennent à $\M_\infty$. Noter que, dans le cas d'un intervalle $E$ de longueur 4 dont les extrémités
 ne sont pas dans $\Z$, on ignore si Irr$_E$ est infini.
 
 \medskip
 
 Le théorème suivant décrit la structure des supports des mesures $\mu \in \M_\infty$ :
 
 \smallskip
     
   \n {\bf Théorème 1.6.4.} {\it Soit $S$ un compact de $K$. Pour que $S$ soit le support d'une mesure
   appartenant à $\M_\infty$, il faut et il suffit que les deux conditions suivantes soient satisfaites} :
   
   (i) \ $\cp(S) \geqslant 1$.
   
   (ii) {\it $S$ est réduit, au sens de} A.4.7.

   \smallskip
   
   \n {\it Démonstration}. La nécessité des conditions (i) et (ii) a déjà été démontrée (th.~1.2.6). Pour la suffisance,
   on applique le th.1.6.3 à $S$; on a $\mu_S \in \M_\infty$, et puisque $S$ est réduit, on a $S=\Supp(\mu_S)$.

\smallskip

\n {\bf 1.6.5. Exemple.}

    Soit $E$ {\it l'ensemble triadique de Cantor} dans $[0,1]$ : ensemble des $\sum_{n=1}^\infty \varepsilon_n 3^{-n}$, avec
    $\varepsilon_n \in \{0,2\}$. Soit $\gamma = \cp(E)$. D'après [Ro 64], §2 et [Ra 95], th.~5.3.7, on a $\gamma \geqslant \frac19$, et d'après
    [RR 07] et [LSN 17], il est très vraisemblable que $ \gamma = 0,22094\ldots$

     {\it L'ensemble $E$ est réduit}. Cela résulte de la prop.~A.4.8; en effet, tout ouvert non vide de $E$ contient un sous-ensemble déduit de $E$ par homothétie par une puissance de $1/3$ suivie d'une translation; un tel ouvert n'est donc pas de capacité $0$. 
     
     Ainsi, on  peut appliquer le th.~1.6.4
    à un multiple  $\lambda E$ de $E$, avec $\lambda > 1/\gamma$, par exemple $\lambda > 9.$ 
        
    L'intérêt de cet exemple est que la mesure $\mu_{\lambda E}$ ainsi obtenue {\it a un support qui est de mesure $0$ pour la mesure de Lebesgue $dx$}. Ainsi, $\mu_{\lambda E}$ et $dx$ sont des mesures {\it étrangères}, au sens de [INT], chap.~V, §5.7:  $\mu_{\lambda E}$ n'est de la forme  $\varphi(x)dx$ pour aucune fonction intégrable $\varphi$ sur $K$. Plus généralement:
    
    \smallskip
     \n {\bf Proposition 1.6.6.} {\it Soit $I$ un intervalle fermé. Pour tout $\varepsilon > 0$, il existe une partie compacte de $I$,
     de mesure de Lebesgue $0$, dont la capacité est $> \cp(I) - \varepsilon.$}
     
    \smallskip
    
    \n {\it Démonstration.} Il suffit de traiter le cas où $I = [-2,2]$; on a alors $\cp(I) = 1$ et $ I$ contient l'ensemble de Cantor $E$
    défini ci-dessus. Soit $n$ un entier $> 0$ et soit $f_n : I \to I$
    l'application $x \mapsto T_n(x)$, où $T_n$ est le $n$-ième polynôme de Chebyshev, cf. (A.2.2). On a $f_n^{-1}(I)=I$. Soit $E_n = f_n^{-1}(E)$. On a $E_n \subset I$, $E_n$ est de mesure de Lebesgue $0$, et, d'après (A.6.1), on a $
    \cp(E_n) = \cp(E)^{1/n}$. Comme on a $0 \ < \cp(E) < 1$, on a donc $\cp(E_n) > 1 - \varepsilon$ pour tout  $n$  assez
    grand. Un tel $E_n$ satisfait aux conditions de la prop.1.6.6.
    
    \smallskip
    
    \n  {\bf Corollaire 1.6.7.} {\it Pour tout intervalle fermé $K$ de longueur $> 4$, il existe $\mu \in M_{\infty,K}$ tel que
    $\Supp(\mu)$ soit de mesure $0$ pour la mesure de Lebesgue.}
    
    \smallskip
    
    \n {\it Démonstration.} Cela résulte de la prop.1.6.6 et du th.1.6.4.
    \medskip
    
    \n {\bf Un problème}
    
    On aimerait avoir davantage de renseignements sur les mesures qui appartiennent à~$\M_\infty$.
    Je me borne à un exemple :
    
     \smallskip 
    (1.6.8) {\it Soit $E$ un intervalle fermé de $\R$, et soit $\nu_E$ sa mesure de Lebesgue (normalisée pour être de masse 1).
    Est-il possible que $\nu_E$ appartienne à} $\M_\infty$ ?
    
      Il me paraît probable que la réponse est «non». En tout cas, une condition nécessaire est que  $I(\nu_E) \geqslant 0$,
    où $I(\nu_E) = \frac{1}{L^2} \iint_{E \times E} \log|x-y| dx dy$, $L$~étant la longueur de $E$. Un calcul élémentaire
    donne $I(\nu_E) = \log L - 3/2$. On doit donc avoir  $L \geqslant e^{3/2} = 4,816...$, ce qui est un peu plus restrictif
    que la borne évidente $L > 4$. Dans certains cas, on peut améliorer cette borne : si par exemple le milieu de $E$ est $0$,
    en écrivant que $\int_E \log |x| dx$ est $ \geqslant 0$, on trouve $L \geqslant 2e = 5,436...$

%
\subsection{Les sous-ensembles d'un cercle centré en $0$} 
    
    On suppose maintenant que $K$ est {\it un cercle $C$ de centre $0$ et de rayon $r > 1$, avec} $r^2 \in \N$. 
    
    \smallskip
    {\small Le cas qui nous intéressera par la suite est celui où $r^2 =  q$, avec $q = |\F|$, comme dans l'introduction. Pour le cas plus général où une puissance entière de $r$ est dans $\N$, voir [Ro 69], §2.\par}
    
    \smallskip 
    Soit $I = [-2r,2r]$. Si $z \in \C$, posons $f(z) = z+\zbar$. Nous obtenons ainsi une application continue surjective
    $f : C \to I$; si $a\in I$, l'image réciproque de $a$ est formée des racines de l'équation $z^2-az+r^2=0$ (on voit ainsi que, si $a$ est un entier algébrique, il en est de même de $z$). Cette application va nous permettre {\it d'identifier les ensembles $\Irr_C, \M_C, \M_{C,\infty}$ avec les ensembles correspondants pour~$I$}. Cela provient des propriétés suivantes de $f$:
        
    \smallskip

    \n 1.7.1. {\it Comportement vis-à-vis des mesures}
    
      \smallskip
      
        L'application $\mu \mapsto f(\mu)$ est un isomorphisme de l'espace des mesures sur $C$ invariantes par conjugaison sur l'espace des mesures sur~$I$.
    
    \smallskip
    
   \n 1.7.2. {\it Comportement vis-à-vis de la capacité }
     
      \smallskip

     Si $E$ est un compact de $C$ stable par conjugaison, soit $E_I$ son image par  $f$. D'après la prop.A.7.1, on a:  
    
     \smallskip

    (1.7.3) \ $\cp(E) = r^{1/2}\cp(E_I)^{1/2}$.
    
     \smallskip

    (1.7.4)  Si $\cp(E) > 0$, l'image par $f$ de la mesure d'équilibre $\mu_E$ est la mesure d'équilibre $\mu_{E_I}$.
    
     \smallskip

      \n  1.7.5. {\it Comportement vis-à-vis des entiers algébriques -- cf.} [Ro 6],~§2
      
        \smallskip
      Si $P \in \Irr_C$ il existe un unique $P_I \in \Irr_I$ dont les racines sont les images par $f$ des racines de $P$, et l'on obtient ainsi une bijection $\Irr_C \to \Irr_I$. On a $\deg P_I = \frac{1}{2}\deg P$, sauf lorsque $r$ est entier et
      $P(X) = X±r$ auquel cas $P_I(X)= X±2r$. De plus, on a :
      
       \smallskip
      (1.7.6) \ $f(\delta_P) = \delta_{P_I}$.

  \smallskip

 \n   Ces propriétés entraînent :
    
      \smallskip
      
      \n {\bf Proposition 1.7.7.} {\it L'application $\mu \mapsto f(\mu)$ définit un isomorphisme de $\M_C$ sur $\M_I$
      qui transforme $\M_{C,\infty}$ en $\M_{I,\infty}$.}

        \smallskip
      En appliquant à $I$ le th.1.6.3, le th.1.6.4 et le cor.1.6.7, on obtient:

        \smallskip
    
    \n {\bf Théorème 1.7.8.} {\it Soit $E$ une partie fermée de $C$ invariante par conjugaison.}
    
     (i) {\it Pour que $E$ soit le support d'une mesure appartenant à $\M_{C,\infty}$, il faut et il suffit  que $E$ soit réduit et $\cp(E) \geqslant r^{1/2}.$}
     
    (ii) {\it Il existe des ensembles $E$ satisfaisant à} (i) {\it  de mesure de Lebesgue $0$.}

     (iii) {\it Si $\cp(E) \geqslant r^{1/2}$, on a $\mu_E \in \M_{C,\infty}$.}
     
     \smallskip
         
    \n{\small [Noter le remplacement de $1$ par $r^{1/2}$, dû à (1.7.3).] }
    
    \smallskip
    
    \n {\it Exemple}.  Soit $L = [a,b]$
      un intervalle contenu dans $[-r,r]$
    et de longueur $\geqslant2$; soit $E_L$ l'ensemble des points $z = x+iy \in C$ tels que $x \in L$. Il résulte de (1.7.3) que $\cp (E_L ) \geqslant r^{1/2}$, de sorte que la mesure d'équilibre $\mu$ de $E_L$ appartient à  $\M_{C,_\infty}$ d'après le th.1.7.8. Sur la « moitié supérieure » de $E_L$ (celle où $y \geqslant 0$), on a $ \mu = \frac{1}{2\pi}\frac{dx}{\sqrt{(b-x)(x-a)}}$; idem pour l'autre moitié.

%
%
\subsection{Application aux variétés abéliennes}

    Comme on l'a dit dans l'introduction, à toute variété abélienne $A \neq 0$ sur un corps~$\F$ à~$q$ éléments,
    on associe le polynôme caractéristique  $P_A$ de son endomorphisme de Frobenius. Ses valeurs propres se trouvent sur le cercle~$C$ de centre~$0$ et de rayon $r=q^{1/2}$. 
    
    \smallskip
  \n {\bf Lemme 1.8.1.} {\it Soit $P \in \Irr_C$. Il existe un entier $m > 0$ et une variété abélienne~$A$ sur~$\F$
  tels que $P_A = P^m.$}
  
     \smallskip
     
     Cela résulte du théorème de Honda-Tate, cf. [Ta 69], §1, Remarque 2.
     
     \smallskip
     Comme $\delta_{P^m} = \delta_P$, on déduit de là :
     
       \smallskip
  \n {\bf Proposition 1.8.2.} {\it L'enveloppe convexe fermée $\M^{\rm ab}$ des $\delta_A$ est égale à $\M_C$.}  
  
     \smallskip
     
     On peut donc appliquer à $\M^{\rm ab}$ les résultats du § précédent, et en particulier le th.1.7.8; noter que $r^{1/2}$ est ici~$q^{1/4}$:
       la capacité du support d'une mesure appartenant à $\M^{\rm ab}_\infty$ est donc $\geqslant q^{1/4}$.
     
     \medskip
     
    
%
%
%
%
\section{Démonstration du théorème de Robinson}

\subsection{Résumé de la démonstration}

\smallskip
Il s'agit de démontrer le th.1.6.2, et donc aussi le th.1.6.1. On se donne des nombres réels  $$a_0 < b_0 < a_1 < \cdots < a_g < b_g,$$ en nombre 
$2g+2, g \geqslant 0$. Pour chaque $j = 0, \ldots, g$, soit $E_j = [a_j,b_j]$ et soit $E = \cup_j E_j$. On doit prouver que, si
$\cp(E) > 1$, il existe une suite $(P_n)$ de polynômes unitaires à coefficients dans $\Z$, dont toutes les racines appartiennent à $E$, et
qui jouissent de la propriété :

(2.1.1) \ $\mu_E = \lim_{n\to \infty} \mu_{_{P_n}}.$

\medskip
La méthode est la suivante:

\smallskip
\n On commence par traiter un cas particulier, celui que nous appellerons {\it de Pell-Abel}, cf. §2.2; dans ce cas, on verra au §2.4 qu'il existe, pour un certain entier $r \geqslant 1$, un polynôme unitaire $P(x)$ de degré  $r$, à coefficients réels, dont les racines sont distinctes et appartiennent à $E$; de plus, sur chaque $E_j$, ce polynôme oscille entre $-M$ et $M$ (avec $M >0$) de façon analogue à celle des polynômes de Chebyshev usuels sur $[-2,2]$. On a $\cp(E)=M/2$, de sorte que $M >2$. 

L'étape suivante
consiste à rétrécir $E$ et à modifier $P$, arbitrairement peu, de telle sorte que l'on ait encore $\cp(E)>1$, et que les coefficients de  $P$  soient rationnels
(§2.6); de là, on parvient à un polynôme à coefficients {\it entiers}, dont les racines sont dans $E$ (§2.7); on obtient ainsi des polynômes de degré arbitrairement grand,  et l'on prouve qu'ils jouissent de la propriété (2.1.1), ce qui démontre le th.~1.6.2 dans le cas considéré. Le cas général se ramène à celui-là en prouvant que, en rétrécissant arbitrairement peu $E$, il devient du type de Pell-Abel (§2.8 et §2.9).

\medskip

  La démonstration résumée ci-dessus est essentiellement celle de Robinson ([Ro 64]). Il en existe une autre, due à R. Rumely, qui est plus longue, mais qui donne un résultat plus général : elle s'applique à des courbes algébriques de genre quelconque et elle permet d'imposer aux points rationnels de ces courbes des conditions, non seulement archimédiennes, mais aussi $p$-adiques. Je renvoie à [Ru 13] pour les énoncés (Introduction, th.0.4 et chap.4, th.4.2), ainsi que pour l'historique du sujet, et notamment le rôle de D. Cantor ([Ca 69], [Ca 80]),  qui a été le premier à donner des énoncés de ce type.
%
%
%
\subsection{Courbes hyperelliptiques et équation de Pell-Abel}

  Soit $k$ un corps de caractéristique $\neq 2$, et soit $D(x) \in k[x]$ un polynôme unitaire de discriminant $\neq 0$, et de degré pair $2g+2, g \geqslant 0$. L'équation $y^2 = D(x)$
  définit une courbe hyperelliptique affine  $C^{\rm aff}_D$ de genre  $g$; sa complétée $C_D$ s'obtient en lui ajoutant deux points à l'infini, que nous noterons $\infty_+$ et $\infty_-$; ils sont caractérisés par le fait que $y/x^{g+1}$ prend la valeur  1  au premier et $-1$ au second.
  
  \smallskip
  L'algèbre affine de $C^{\rm aff}_D$ est un $k[x]$-module libre de base $\{1,y\}$. Ses éléments inversibles
  sont de la forme:
  
    \smallskip
  (2.2.1)  $f = P + yQ$, \ avec $P,Q \in k[x]$  \ et $P^2 - DQ^2 = c$, \ avec $c \in k^\times$.
  
    \smallskip
 On supposera que $P$ et $Q$ sont non nuls (cela revient à éliminer le cas $f \in k^\times$), et on les normalisera en demandant que ce soient des polynômes unitaires. On appellera {\it degré} de $f$ le degré de $P$; on a deg$(P) \geqslant g+1$. L'équation  $P^2 - DQ^2 = c$ est l'analogue pour $k[x]$ de l'équation dite «de Pell»
 pour $\Z$. Il est raisonnable de l'appeler {\it l'équation de Pell-Abel}, car elle apparaît pour la première fois dans Abel ([Ab 26]) où elle est étudiée à l'aide de fractions continues, comme l'équation de Pell dans le cas de~$\Z$.
  
 On trouve aussi dans la littérature le nom d'{\it équation d'Abel-Chebyshev}: elle avait en effet été retrouvée, et utilisée, par
 Chebyshev ([Ch 54]), à propos d'un problème de mécanique\footnote{Problème (important pour les constructeurs de locomotives): comment utiliser certains quadrangles articulés
 (les {\it mécanismes de Chebyshev}) pour transformer aussi bien que possible un mouvement circulaire en un mouvement rectiligne, et inversement ? C'est en essayant d'optimiser le «aussi bien que possible» que Chebyshev a été conduit aux polynômes qui portent son nom, ainsi qu'à l'équation $P(x)^2 - D(x)Q(x)^2 = c$. Le lecteur curieux trouvera
 sur internet des reproductions (avec vidéo) de certains de ces mécanismes.}.

\bigskip
 Une différence importante avec l'équation de Pell usuelle est que {\it l'équation de Pell-Abel n'a pas toujours de solution}.
 
 \bigskip
 De façon plus précise, soit $r$ un entier $\geqslant 1$. Il y a équivalence entre :
 
 (2.2.2) {\it L'équation (2.2.1) a une solution de degré $r$.}
 
 (2.2.3) {\it Le diviseur $r(\infty_- - \infty_+)$ de la courbe $C_D$ est linéairement équivalent à}  $0$.
 
 \smallskip
 {\small [C'est immédiat : il suffit de vérifier qu'une fonction rationnelle $f$ sur $C_D$ de diviseur $r(\infty_- - \infty_+)$
 est nécessairement de la forme  $f = P + yQ$, comme dans (2.2.1).]\par}
 
 \smallskip
  On peut reformuler (2.2.3) de la manière suivante: si $J$ désigne la jacobienne de $C_D$, le diviseur $\infty_- - \infty_+$ définit un
  point $P_\infty$ de $J(k)$, et (2.2.3) signifie que  $rP_\infty = 0$, autrement dit que $P_\infty$ est un point d'ordre fini divisant $r$. Lorsque $g=0$, c'est le cas car $J=0$, et l'on trouve pour  $P,Q$  les polynômes de Chebyshev usuels (de première et seconde espèce, respectivement), avec une normalisation un peu différente;
 c'est aussi le cas (pour $r$ bien choisi) lorsque $g \geqslant 1$ et $k$ est un corps fini. Par contre, si $g\geqslant 1$ et si $\car(k)=0$,
 il est facile de construire des exemples où $P_\infty$ est d'ordre infini, ce qui entraîne que l'équation de Pell-Abel n'a pas de solution.
 
 \medskip
 
 \n {\small [Une autre façon de formuler ceci est d'introduire la {\it Jacobienne généralisée}  $J_{\mathfrak m}$ de $C_D$ relative au conducteur $\mathfrak{m} = \infty_- \ +  \infty_+.$ On a une suite exacte $ 1 \to \G_m \to J_{\mathfrak m} \to J \to 1$,
 
 \n qui montre que $J_{\mathfrak m}$ est une extension de $J$ par le groupe multiplicatif $\G_m$; la classe de cette extension  dans $\Ext(J,\G_m) \simeq J(k)$
 est $P_\infty$, au signe près. Dire que $P_\infty$ est d'ordre fini signifie donc que $J_{\mathfrak m}$ {\it est isogène à}  $J \times \G_m$, ou encore que le 1-motif mixte de la courbe $C^{\rm aff}_D$ est scindé.]\par}

%
%
%
\subsection{Le cas réel: la forme de troisième espèce canonique}

Nous allons maintenant supposer que $k = \R$ et que le polynôme $D(x)$ a pour racines les $a_j,b_j$ du §2.1 :

\smallskip
(2.3.1) \hspace{5mm}  $a_0 < b_0 < a_1 < \cdots < a_g < b_g.$

\smallskip
Ainsi $D(x)$ est $\leqslant 0$ si $x \in E$ et $> 0$ sinon. Les points réels de la courbe $C_D$ correspondent donc aux
$x \in \P_1(\R) \sm \overset {\;\circ} E$.

\smallskip

  Soit $A \in \R[x]$ un polynôme unitaire de degré $g$; on lui associe la forme différentielle $\eta_A = \frac{A(x)dx}{y}$
  sur $C_D$. Cette forme est holomorphe ailleurs qu'en l'infini; elle a un pôle simple en $\infty_+$ et $\infty_-$, avec résidus $-1$ et $+1$ respectivement; c'est une «forme de 3-ème espèce». Changer $A$ revient à lui ajouter une forme de première espèce, puisque celles-ci ont pour base les $x^jdx\! / \!y, 0\leqslant j < g$. On en déduit {\it qu'il existe un choix
  et un seul de  $A$ tel que les périodes réelles de $\eta_A$ soient nulles}, autrement dit qui vérifie les conditions :
  
  \medskip
   (2.3.2)   \  $\int_{b_{j-1}}^{a_j} \frac{A(x)}{\sqrt{D(x)}}\,dx = 0$ \ pour $j=1,\ldots,g.$

\smallskip

\n {\small[Dans cette formule, ainsi que dans les suivantes, si  $t$ est réel $\geqslant 0$, nous notons  $\sqrt{t}$ sa racine carrée $\geqslant 0$; si $t<0$, nous définissons $\sqrt{t}$ comme $i\sqrt{-t}$.]\par}

\smallskip

  Le polynôme $A$ déterminé par ces conditions sera noté $R$, et la forme  $\eta_R $ sera appelée la
  {\it forme de 3-ème espèce canonique}, et notée simplement $\eta$. La formule (2.3.2), appliquée à $R$, montre que,
  dans chaque intervalle intermédiaire $T_j = \ [b_{j-1},a_j]$, l'intégrale de $\frac{R(x)}{\sqrt{D(x)}}$ est nulle; cela entraîne que $R(x)$ change de signe dans $T_j$, donc s'annule en au moins un point intérieur de $T_j$ . Comme le nombre des $j$  est égal au degré de $R$, on obtient :
  
  \smallskip
  (2.3.4) {\it Le polynôme $R(x)$ a une racine et une seule à l'intérieur de chaque $T_j$, et n'a aucune autre racine $($réelle ou complexe$)$.} 
  
  \smallskip
  Les changements de signe de  $R$  se font donc dans les «trous» $T_j$. Cela montre que $R$ garde un signe constant sur chaque $E_j$. Comme $R(x)$ est $\geqslant 0$ pour $x$  assez grand, on en déduit, par récurrence descendante sur $j$:
  
  \smallskip
  
  (2.3.5) {\it On a $R(x) \geqslant 0$  sur $E_j$ si $g-j$ est pair et  $R(x) \leqslant 0$ sinon.         }

  \medskip
  
  On va  maintenant s'intéresser aux «périodes imaginaires» de $\eta$, qui sont égales à~$2\eta_j$, avec :
  
    \medskip
  
  (2.3.6)  $\eta_j = \int_{a_j}^{b_j} \frac{R(x)}{\sqrt{D(x)}} \,dx$, $j=0,1,\ldots,g . $
  
  \medskip
  
  Leurs parties réelles sont nulles.
  
  \bigskip
  
  \n {\bf Proposition 2.3.7.} {\it Il existe des signes $\varepsilon_j \in \{-1,1\}$
  tels que l'on ait :}
  
  \smallskip
   (2.3.8)  $\sum_{0\leqslant j \leqslant g}  \varepsilon_j \eta_j = i\pi$.
   
\medskip

\n {\it Démonstration.} Notons $S$ la variété analytique complexe $C_D(\C)$, vue comme une surface de Riemann au-dessus de $ \P_1(\C) = \C \cup \{\infty\}$ via l'application $(x,y) \mapsto x$. C'est un revêtement quadratique, ramifié seulement aux points $a_j,b_j$. Ce revêtement est non ramifié (c'est un vrai revêtement, au sens de la Topologie)
au-dessus de $U = \P_1(\C) \sm E$. Notons  $\tilde{E_j}$  l'image réciproque de $E_j$ dans $S$, et soit $\tilde{E}$ la réunion des $\tilde{E_j}$. L'image réciproque $\tilde{U}$ de $U$ dans $S$ est égale à $S \sm \tilde{E}$. Un calcul sur les groupes fondamentaux montre que {\it le revêtement $\tilde{U} \to U$ est trivial}.

\smallskip
{\small [En termes plus traditionnels, cela signifie que l'on peut définir $y$ comme une fonction méromorphe sur $U$, par prolongement analytique à partir de l'une de ses deux valeurs possibles en un point de $U$.]\par}

\smallskip
Il en résulte que $\tilde{U}$ a deux composantes connexes; nous noterons $\tilde{U}_+$ celle qui contient le point $\infty_+$; définition analogue pour $\tilde{U}_-$. L'adhérence $S_+$ de  $\tilde{U}_+$ est $\tilde{U}_+ \cup \tilde{E}$; c'est une surface à bord (une «pièce», dit Bourbaki [FRV], 11.1.2), dont le bord est la réunion des $g+1$ cercles $\tilde{E}_j$.
Si l'on définit de manière analogue $S_-$, on voit que $S$ peut s'obtenir à partir de $S_+$ et $S_-$ en collant leurs bords. 

On notera que l'orientation naturelle de $\tilde{U}_+$ (comme variété analytique complexe) définit une orientation
sur chacune des composantes $\tilde{E}_{j+}$ de son bord; cela donne un sens à une intégrale du type $\int_{\tilde{E}_{j+}} \omega$. Lorsque $\omega$ est de première espèce, la formule de Cauchy\footnote {Cas particulier de la formule de Stokes en dimension 2.} ([FRV], 11.2.5) donne :

\smallskip
  $ \sum_j \int_{\tilde{E}_{j+}} \omega = 0. $
 
 \smallskip
 Un argument analogue, appliqué à $\eta$ et à $\tilde{U}_+$ dont on a retiré un petit disque autour de $\infty_+$,
 donne :
 
 \smallskip
 (2.3.9) $ \sum_j \int_{\tilde{E}_{j+}} \eta= -2i\pi$,
 
 \smallskip
 
\n du fait que le résidu de $\eta$ au point $\infty_+$ est~$-1$.

\smallskip

On passe de (2.3.9) à (2.3.8) en remarquant que, pour tout $j$, l'intégrale $\int_{\tilde{E}_{j+}} \eta$ est égale à $±2\eta_j$, le signe « $±$ » provenant du fait que l'on n'essaie pas de préciser les orientations.

\medskip
%
%
%
%
\subsection{L'équation de Pell-Abel dans le cas du §2.3 -- énoncé des résultats}

\medskip
On conserve les hypothèses du §2.3, et l'on suppose en outre que, pour un certain entier $r \geqslant 1$,
il existe une solution $(P,Q)$ de degré~$r$ de l'équation  $P^2 - DQ^2 = c$, avec $c \in \R^\times$.
La constante $c$ est positive, car c'est le carré de la valeur de $P$ en l'un quelconque
des points $a_j,b_j$. Nous la noterons $M^2$, avec $M > 0$. On a donc :

\bigskip

  (2.4.1) $P^2 - DQ^2 = M^2$.
  
  \medskip
  
  Supposons $x$ réel; on a $D(x) \leqslant 0$ si et seulement si $x \in E$; cela entraîne :
  
  \medskip
  (2.4.2) $|P(x)| \leqslant M \ \Longleftrightarrow \ x \in E$ ou $Q(x)=0.$
  
  \bigskip
  
  \n [En fait, le th.~ 2.4.4 ci-dessous entraîne que $Q(x) \neq 0$ si $x \notin E$. On peut donc supprimer «ou $Q(x)=0$» dans (2.4.2).]

 \medskip 
 
 De plus : 
  
  \medskip
  (2.4.3) $|P(x)| = M  \Longleftrightarrow \ x$ est, soit une racine de $Q$, soit l'un des $a_j,b_j$.

\bigskip

En ce qui concerne les racines des polynômes $P, Q$, la propriété la plus importante pour la suite est la partie (i) du théorème suivant :

\bigskip

\n {\bf Théorème 2.4.4.} (i) {\it Les racines des polynômes $P$ et $Q$ sont réelles, simples, et appartiennent à l'intérieur $\overset {\;\circ} E$ de $E$.}

(ii) {\it Soit $f = P + yQ$, cf. §2.2.1. On a $df\! / \!f = r\eta$, où $\eta$ est la forme de troisième espèce canonique définie au §2.3.}

\bigskip

   En fait, on aura besoin d'autres résultats, à savoir :
   
   \bigskip
   
   \n {\bf Théorème 2.4.5.} {\it  Pour $j =0,\ldots,g$, soit $r_j = {r\over {\pi}}|\eta_j|$, où $\eta_j$ est la $j$-ième demi-période imaginaire
   de $\eta$, cf. (2.3.6). Alors }:
      
      (2.4.6) {\it Les $r_j$ sont des entiers $> 0$, de somme $r$.}
      
      (2.4.7) {\it Le nombre de racines de $P$ dans $E_j$ est $r_j;$ le nombre de racines de $Q$ dans $E_j$ est $r_j-1$.}
      
      (2.4.8) {\it Les racines de $Q$ qui sont contenues dans $E_j$ divisent $E_j$ en $r_j$ sous-intervalles$;$ dans chacun d'eux,
      le polynôme $P$ est, soit strictement croissant, soit strictement décroissant, de valeurs extrêmes  $M$
      et $-M$.}

\medskip

 \ \ La figure ci-après indique à quoi ressemble le graphe de $P$ dans le cas où $g=1,\\ r_1=~4, r_2=6, r=10$, les racines de $Q$ étant notées $c_1,\ldots, c_8$. Les dix sous-intervalles sont $[a_0,c_1], [c_1,c_2], \ldots, [c_8,b_1]$.
 
 \medskip

\begin{center}

\includegraphics[keepaspectratio,viewport=70 105 450 345,clip,width=\textwidth]{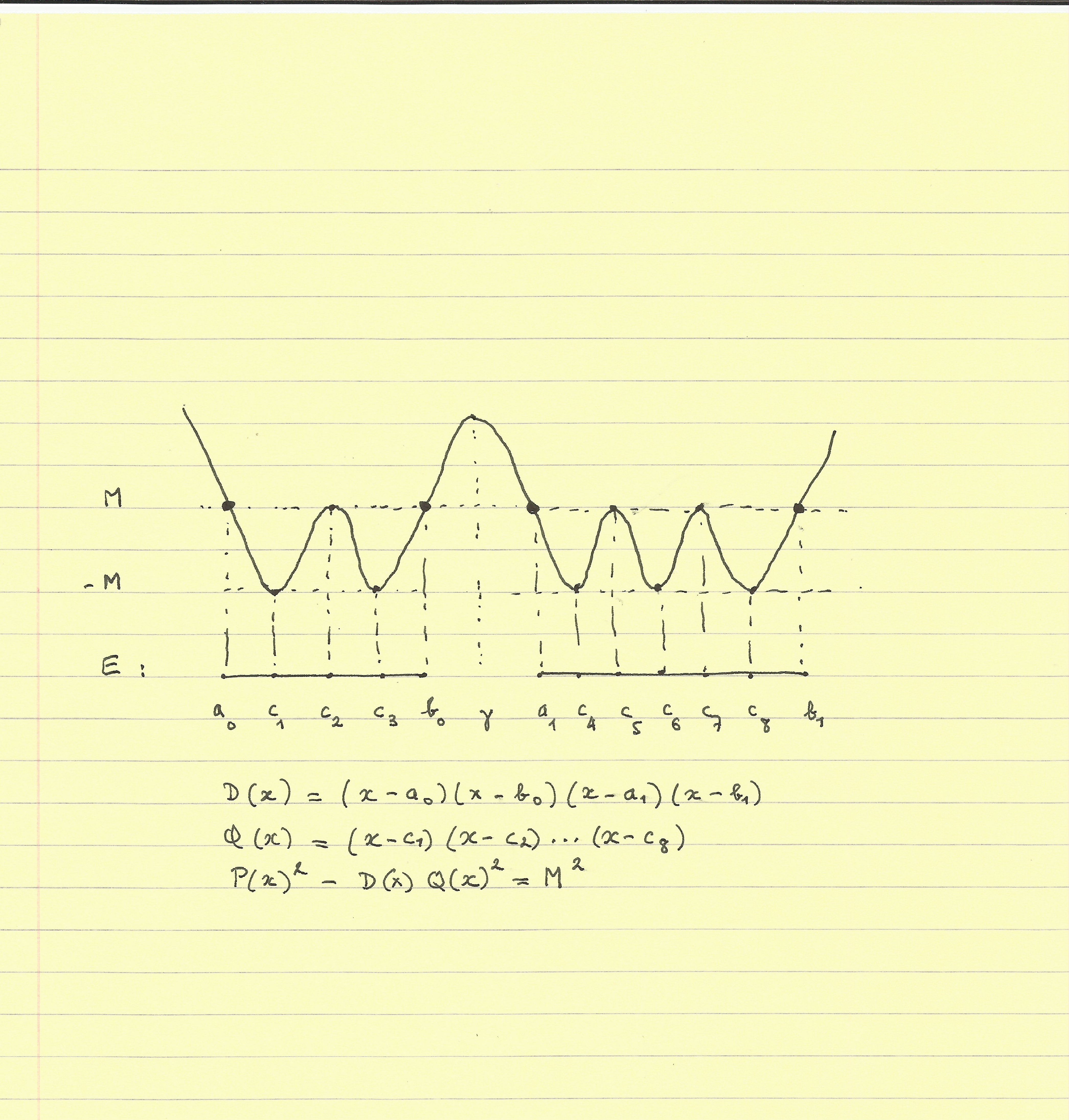}

\end{center}

\medskip

  Autres propriétés :
  
  \medskip

\n{\bf Théorème 2.4.9.} (i) $ dP\! / \!dx = rQ(x)R(x)$, {\it où $R$ est le polynôme de degré $g$ défini au~§2.3.}
 
 (ii) {\it Le polynôme $P$ est le $r$-ième polynôme de Chebyshev {\rm(au sens du §A.2)} du compact~$E$.}
 
 (iii) $\cp(E) = (M/2)^{1/r}.$

\smallskip
\n [Dans le cas de la figure ci-dessus, on a $R(x) = x - \gamma$, où $\gamma$ est le point de $[b_0,a_1]$ où $P$ est maximum.]

\medskip

Les démonstrations des théorèmes  2.4.4, 2.4.5 et 2.4.9 seront données au §2.5.

\smallskip
{\small {\it Remarque historique.} Il n'est pas facile de dire à qui ces théorèmes sont dus. 
 L'une des difficultés est qu'ils sont rarement énoncés explicitement,
et, du coup, ils ne sont pas démontrés en détail. La plus ancienne référence que
j'aie pu trouver est celle de Robinson [Ro 64]. Il y a eu ensuite  Peherstorfer ([Pe 90]), Sodin-Yuditskii ([SY 92]), Bogatyrëv ([Bo 99] et [Bo 05]), et sans doute d'autres. Il est toutefois possible
que l'école russe, descendante de Chebyshev et Zolotarëv, et en particulier N.I. Akhiezer, ait été familière avec ces résultats dès les années 1930.\par}

%
%
\subsection{Démonstrations des énoncés du §2.4}

 \n {\it Démonstration du th.~2.4.4 (ii).} Soit $f = P +yQ$; les seuls pôles de $df\! / \!f$ sont les points
  $\infty_+$ et $\infty_-$, qui sont des pôles simples de résidus $-r$ et $r$ respectivement. Pour prouver que $df\! / \!f=r\eta$, il nous suffit donc de prouver que les périodes réelles de $df\! / \!f$ sont nulles. 
  
  \smallskip
 Si $j = 1,\ldots,g$, choisissons une détermination de $y$ sur l'intervalle $T_j= [b_{j-1},a_j]$ , par exemple $y(x) = \sqrt{D(x)}$; cela permet de voir
  $f$ comme une fonction $f(x)$ sur $T_j$. Cette fonction est réelle, et ne s'annule pas; comme $|f(b_{j-1})| = M = |f(a_j)|$, on en conclut que $f(b_{j-1}) = f(a_j) = \varepsilon M$, avec $\varepsilon = ±1$; le signe de $f$ est donc $\varepsilon$ sur $T_j$  La fonction $\log(\varepsilon f)$ prend les mêmes valeurs aux extrémités de $T_j$. Cela entraîne que l'intégrale sur  $T_j$ de sa dérivée est $0$;
  comme cette dérivée est $df\! / \!f$, cela signifie que la $j$-ième période réelle de $df\! / \!f$ est nulle. 
  
   \smallskip
\n {\it Remarque.}  Soit $\varphi(x) = \log(\varepsilon f(x)/M)$. On a
$f(x) = \varepsilon Me^{\varphi(x)}$;  si $f_1 = P -yQ$, le fait que $f_1f = M^2$ montre que $f_1(x) = \varepsilon Me^{-\varphi(x)}$. 
Comme $P = (f + f_1)/2$ et $yQ = (f-f_1)/2$, cela donne :

\smallskip

\noindent (2.5.1)
{\it On a $P(x) = \varepsilon M\cosh(\varphi(x))$ et  \ $y(x)Q(x) = \varepsilon M\sinh(\varphi(x))$ si} \mbox{$x \in T_j=  [b_{j-1},a_j]$.}

\smallskip
{\small Comme $\varphi$ est une «primitive» de $df\! / \!f = rR(x)dx\! / \!y(x)$, on obtient donc une expression explicite de $P$ et $Q$ en termes de $R(x)$, lequel peut se calculer en résolvant un système linéaire d'équations 
de type $g \times g$.\par}

\medskip

\n {\it Démonstration du th.~2.4.9 (i).} 

\smallskip

On a $P + yQ = f$ et $P-yQ ={M^2 \over f}$, d'où  $P= {1 \over 2}(f+{M^2 \over f}),\  yQ={1 \over 2}(f-{M^2 \over f})$,\\
et \ $dP = {1 \over 2}(1- {M^2 \over f^2} )df =yQdf\!/\!f =yQr\eta=rQRdx.$   \medskip

\n {\it Démonstration du th.~2.4.5 et du th.~2.4.4 (i).} 

\smallskip

\n {\small [La méthode est essentiellement la même que pour le th.~2.4.4 (ii), à cela près que le groupe multiplicatif $\R^\times$ est remplacé par le cercle unité $C_1$.]\par}

\smallskip
Considérons l'un des intervalles $E_j$, $j = 0,\ldots,g$. Posons $y_1(x) =
\sqrt{-D(x)}$, qui est à valeurs réelles $\geqslant 0$, et choisissons comme détermination de  $y$ la fonction $iy_1(x)$.
Comme ci-dessus, on peut voir  $f(x) = P(x) + iy_1(x)Q(x)$ comme une fonction sur $E_j$; cette fonction est
à valeurs complexes et de module $M$. On peut donc l'écrire sous la forme $f(x) = Me^{i\vartheta(x)}$, où $\vartheta(x)$
est bien défini mod \ $2\pi$. On peut relever $\vartheta: E_j \to \R/2\pi \Z$ en une fonction continue $\theta: E_j \to \R$,
qui est bien déterminée si on lui impose sa valeur en~$a_j$; comme $f(a_j)=± M$, cette valeur est de la forme
$c_0\pi$, avec $c_0 \in \Z$; la valeur de~$f$ en~$b_j$ est alors $c_1\pi$ avec $c_1\in \Z$. On a :

\medskip
(2.5.2) \quad $P(x) = M\cos(\theta(x))$ \ et \ $y_1(x)Q(x) = M\sin(\theta(x)).$

\medskip

\n De plus :

\medskip
\quad $ \int_{a_j}^{b_j}  df\!/\!f = i  \int_{a_j}^{b_j} d\theta = i(\theta(b_j)- \theta(a_j)) =(c_1 - c_0)i\pi$. 

\medskip

Soit $r_j = |c_1-c_0|$; c'est un entier $>0$.
\n Comme $df\!/\!f = r\eta$, la formule ci-dessus montre que, au signe près, $r_ji\pi$ est la $j$-ième demi-période imaginaire
de $r\eta$. On a donc $r_j = r|\eta_j|/\pi$, avec les notations de (2.3.3).  D'après (2.3.5), il existe des signes~$±$
tels que :

\smallskip

(2.5.3)  \quad  $ \sum_j  ± r_j = r$.

\smallskip

Noter que la dérivée de $\theta(x)$ ne s'annule pas, donc  $\theta(x)$ est strictement croissante (resp. décroissante)
si $c_0 < c_1$ (resp. si $c_1 < c_0$). Or, si une variable croît strictement entre deux multiples entiers $u\pi$ et $v\pi$, son cosinus s'annule un nombre de fois égal à $v-u$, et son sinus s'annule $v-u-1$ fois en dehors des extrémités. On en conclut
que {\it le nombre de racines de $P$ dans $E_j$ est $r_j$ et que le nombre analogue pour $Q$ dans l'intérieur de $E_j$ est $r_j-1$}. Comme $P$ a au plus $r$ racines, on a donc :

\smallskip
(2.5.4)  \quad  $\sum_j r_j \leqslant r.$

\smallskip

En comparant (2.5.3) et (2.5.4), on voit que tous les signes de (2.5.3) sont des signes~$+$, de sorte que l'on a
$\sum r_j = r$, et l'on voit aussi que $P$ n'a pas d'autres racines que celles qui sont dans les $E_j$ et que celles-ci sont des racines simples (ce dernier point résulte aussi de la formule $dP/dx = rQR$ du th.~2.4.9 (i)). Le même argument montre que $Q$ n'a pas d'autres racines que celles contenues dans l'intérieur de $E$, et que celles-ci sont simples. Cela achève la démonstration de 2.4.4 (i), (2.4.6) et (2.4.7).
Pour (2.4.8), on remarque que, dans un sous-intervalle $[\alpha, \beta]$ du type de (2.4.8), la dérivée de $P$ est
$\neq 0$ en tout point $x$ tel que $\alpha < x < \beta$: cela résulte de la formule $\frac{dP}{dx}= rQR$ et du fait que $R$
ne s'annule pas sur les $E_j$, cf. (2.3.4).

\medskip

\n {\it Démonstration du th.~2.4.9 (ii) $:$ $P$ est un polynôme de Chebyshev.} 

\nobreak\smallskip\nobreak

On démontre d'abord:

\smallskip

\n {\bf Lemme 2.5.5.} {\it Soit $q \in \R[x]$ tel que $|q(x)| < M$ pour tout $x\in E$. Le polynôme $P-q$ a au moins $r$ racines dans} $E$.

\smallskip 
  
  \n {\it Démonstration.} Soit  $[\alpha, \beta]$ un sous-intervalle de type (2.4.8); on a $P(\alpha) = - P(\beta) =± M$;
  l'hypothèse faite sur $q$ entraîne que le polynôme $P-q$ a le même signe que $P$ en~$\alpha$ et en~$\beta$;
  il a donc au moins une racine dans l'intervalle ouvert $\mathopen]\alpha,\beta\mathclose[$; d'où le lemme, puisque le nombre des sous-intervalles est $r$.
  
  \smallskip

\n {\it Démonstration du th.~2.4.9 (ii).}

Soit $q$ le polynôme de Chebyshev de $E$ de degré  $r$; le fait qu'il soit unique montre que ses coefficients sont réels. Si l'on avait $q \neq P$, on aurait  $\sup_{x\in E} |q(x)| < M.$ D'après le lemme ci-dessus,  le polynôme $P-q$ a au moins $r$ racines; comme il est de degré $<r,$ ce n'est possible que s'il est nul, ce qui contredit l'hypothèse faite sur $q$.

 \medskip
 
 \n {\it Démonstration du th.~2.4.9 (iii) $:$ calcul de $\cp(E)$.} 

\smallskip
On a vu ci-dessus que, dans tout intervalle $E_j$, le polynôme $P$ a $r_j$ zéros. Le même argument montre, pour tout $z \in [-M,M]$, distinct de $±M$, l'équation $P(x) = z$ a $r_j$ solutions dans $E_j$, donc n'a aucune solution
dans $\C \sm E$ puisque $r = \sum r_j$; lorsque $z = ±M$, cet énoncé est également vrai (il résulte par exemple de $P^2-M^2 = DQ^2$). D'où:

\smallskip

(2.5.6) {\it L'image réciproque de $[-M,M]$ par $P : \C \to \C$ est $E$}.

\smallskip

Comme la capacité de $[-M,M]$ est $M/2$, la formule (A.5.7) montre que $\cp(E) = (M/2)^{1/r}.$

\medskip

\n {\small {\it Autre démonstration}: on verra au §2.7 que, pour tout entier $n > 0$, l'équation de Pell-Abel pour $D$, de degré  $rn,$ a une solution $P_n,Q_n,M_n$ avec  $M_n = 2^{1-n}M^n$. Avec les notations de (A.2), cela entraîne que
$c_{rn}(E) = M_n$. D'après (A.2.1), on a :

\smallskip
\n \ \  $\cp(E) = \lim_{n \to \infty} c_{rn}(E)^{1/rn} =  \lim_n (M_n)^{1/rn} = M^{1/r}\lim_n 2^{(1-n)/rn}= (M/2)^{1/r}$.\par }

%
%
\subsection{Comment remplacer les coefficients réels par des coefficients rationnels}

  On conserve les notations $E,D,P,Q,r,M$ des §§ ci-dessus. On désire démontrer :
  
  \smallskip

\n{\bf Proposition 2.6.1.} {\it Pour tout $M'$ tel que $0 < M' < M$ il existe $E'$ et $P'$ \footnote{On aura soin de ne pas confondre  $P'$  et $D'$ avec les dérivées $dP/dx$ et $dD/dx$.}tels que le sextuplet $E',D',P',1,r,M'$ ait les mêmes
propriétés que $E,D,P,Q,r,M$, et en outre} :
  
  \smallskip
  (2.6.2) {\it On a $P' \in \Q[x]$ et} $E' \subset  \overset {\;\circ} E$.

\smallskip

\n [Noter que, dans cette modification, le polynôme $Q$ a été remplacé par 1.] 

\smallskip
\n {\it Démonstration}. 

Soient $u_1 < u_2 < \cdots < u_{2r}$ les racines de $P(x)^2=M'^2$. D'après (2.5.5), les  $u_i$ sont distincts et appartiennent à l'intérieur de $E$; de façon plus précise, pour tout $j$ impair, l'intervalle $[u_j,u_{j+1}]$ est contenu dans l'intérieur du $j$-ième sous-intervalle de $E$ au sens de (2.4.8). Si $P'$ est assez voisin de $P$
(au sens de la topologie naturelle de l'espace des polynômes unitaires de degré $r$), la même propriété vaut
pour les racines de $u'_1,\cdots,u'_{2r}$  de $P'-M'$ et $P'+M'$, et les relations d'ordre entre racines de $P'-M'$ et racines
de $P'+M'$ sont les mêmes que pour $P$. Choisissons un tel $P'$ qui soit à coefficients dans $\Q$. Posons $Q' = 1$ et $D' = P'^2-M'^2$; l'équation de Pell-Abel  $P'^2 - D'Q'^2 = M'^2$ est satisfaite. On a $D' = 
\prod (x-u'_i)$. L'ensemble $E'$ des points où $D' $ est $\leqslant 0$ est la réunion des $[u_j,u_{j+1}]$, $j$ impair.
Il est contenu dans $ \overset {\;\circ} E$: les conditions (2.6.2) sont donc satisfaites.

\smallskip

\n{\bf Corollaire 2.6.3.} {\it On peut choisir $E'$ et $P'$ de telle sorte que $M'$ soit rationnel et que $\cp(E) - \cp(E')$ soit aussi petit que l'on veut.}

\smallskip

C'est clair, puisque  $\cp(E') = (M'/2)^{1/r}$ d'après le th.~2.4.9 (iii) appliqué à $E'$.

\smallskip

\n {\it Remarque.} Le compact $E'$ a $r$ composantes connexes, et chacune d'elles contient une racine de $P'$ et une seule; le genre de la courbe hyperelliptique correspondante est en général $> g$.

%
%
\subsection{Comment remplacer les coefficients rationnels par des coefficients \mbox{entiers}}

  \n {\bf Théorème 2.7.1.} {\it Soient $E,D,Q,P,r,M$ comme dans les §§2.3 et 2.4. Supposons que $M$ et les coefficients de $P$
  sont rationnels, que $Q=1$ et que $M > 2$} (autrement dit $\cp(E) > 1$). {\it
   Il existe une suite $(P'_n)$ de polynômes unitaires, à coefficients dans $\Z$, de degrés
  tendant vers l'infini, ayant les propriétés suivantes} :
  
   (i) {\it Toutes leurs racines sont simples et appartiennent à $E$.}
  
  (ii) {\it On a $\mu_E = \lim \mu_{P'_n}$.}
  
  \medskip
  
  \n {\bf Corollaire 2.7.2.} {\it Les théorèmes 1.6.1 et 1.6.2 sont vrais pour  $E$.}

 \smallskip
 
 \n {\it Démonstration de} (i) (d'après [Ro 64], §6).
 
 \smallskip
  Le point de départ consiste à utiliser la solution donnée $(P,Q,M)$ de l'équation de Pell-Abel pour en fabriquer d'autres de degré $nr$ pour tout $n \geqslant 1$ : il suffit d'élever $P + yQ$ à la puissance $n$, et de regrouper
  les termes en tenant compte de l'équation  $y^2=D$. Posons $\lambda = M/2$, de sorte que $\lambda > 1$. On trouve :
  
  \smallskip
  $(P + yQ)^n = 2^{n-1}(P_n + yQ_n)$ \ et \  $P_n^2-DQ_n^2 = 2\lambda^n$, 
  
  \smallskip
  
  \n où $P_n$ et $Q_n$ sont unitaires de degrés $nr$ et $nr-g-1$, respectivement.
  
  \smallskip
\n   Pour $n=2$, cela donne:
  
  \smallskip
   $(P+yQ)^2 = 2(P_2 + yQ_2)$, avec $P_2 = P^2 - 2\lambda^2$ et $Q_2 = PQ$.
  
  \smallskip
  
\n   Pour $n$ arbitraire, la formule analogue est :
  
 \smallskip  
  (2.7.3)  \ $P_n(x) = \lambda^nT_n(P(x)/\lambda)$,
   
   \smallskip
  \n où $T_n$ est le $n$-ième polynôme de Chebyshev, cf. (A.2.2); cela se démontre, par exemple, en se plaçant dans un sous-intervalle, et en remarquant que, d'après (2.5.2), on a $P(x) = 2\lambda \cos(\theta(x))$ et $P_n(x) = 2\lambda^n \cos(n\theta(x))$. 
  
  \smallskip
  Nous aurons besoin d'une formule pour $T_n$ qui mette en évidence les propriétés de divisibilité de ses coefficients.
  Robinson donne la suivante (démontrée dans [Ro 62], §2):
  
  \medskip
  
   $ (2.7.4) \ T_n(X) = X^n + \sum_{k=1}^{[n/2]} (-1)^k \frac{n}{k}{n-k-1\choose k-1}X^{n-2k}.$
  
  \smallskip
\n   D'où :
  
  \smallskip
  
  (2.7.5) \ $P_n(x) = P(x)^n + \sum_{k=1}^{[n/2]} (-1)^k \frac{n}{k}{n-k-1\choose k-1}\lambda^{2k}P(x)^{n-2k} .$
  
  \medskip
  
\n  Cela peut s'écrire en abrégé :
  
  \smallskip
  (2.7.6) \ $P_n(x) = x^{nr} + \sum_{k >0} \alpha_k x^{nr-k}$,
  
  \smallskip
  
 \n  où, pour chaque  $k > 0, \alpha_k$ est un polynôme à coefficients dans $\Q$ en $n$, en $\lambda$, et en les
  coefficients de $P$; comme $\lambda$ et les coefficients de $P$ sont rationnels, on voit que $\alpha_k$ {\it est   
  un polynôme en $n$, à coefficients dans $\Q$, et de terme constant $0$}.
  
   \smallskip

\n  Posons 
  
   \smallskip

  (2.7.7) \ $A = \sup_{x\in E}(1 + |x| + \cdots + |x|^{r-1}).$
  
   \smallskip

\n  Choisissons un entier $\l > 0$ tel que :
  
   \smallskip

  (2.7.8) \ $\lambda^\l(\lambda -1) \geqslant A/2.$
  
   \smallskip
   
   \n C'est possible puisque $\lambda > 1.$
   
   \smallskip
 
 \n Choisissons un entier $m \geqslant 1$ tel que, pour $k=1,\ldots,\l r$, les $\alpha_k$, considérés comme polynômes en $n$, soient à coefficients dans $\frac{1}{m}\Z$. 
  
   \medskip
\n {\bf Lemme 2.7.9.} {\it Pour tout $n > 0$ divisible par $m$, il existe $C_n \in \R[x]$, de degré $<nr$, tel que} :

(i) \ $|C_n(x)| < 2\lambda^n$ {\it pour tout} $x\in E$. 
  
(ii) \ {\it Les coefficients de $P_n-C_n$ sont des entiers.}
  
   \smallskip

Ce lemme entraîne la partie (i) du th.~2.7.1 : en effet, le polynôme $P'_n = P_n-C_n$ est à coefficients entiers, et,
d'après le lemme 2.5.5, appliqué à $P_n$ et $C_n$, il a $nr$ racines distinctes contenues dans $E$.

\smallskip

\n {\it Démonstration du lemme 2.7.9.}

 On remarque d'abord que les $\l r$ premiers coefficients $\alpha_1, \ldots, \alpha_{\l r}$ de $P_n$ sont entiers; cela provient du fait que $n$ est un multiple de $m$. On peut donc se borner à prendre pour $C_n$ un polynôme de degré $\leqslant nr - \l r-1$. On utilise pour cela la base formée\footnote{Il ne  faut surtout pas utiliser la base $1,x,x^2, \ldots, x^{nr - \l r-1}$,  qui conduirait à une majoration inutilisable de $|q(x)|$.} des $x^jP_k(x)$,
avec $0\leqslant j < r$ et $0\leqslant k < n-\l$. On écrit $C_n$ sous la forme  $C_n(x) = \sum_{j,k} c_{jk}\ x^jP_k(x)$, et l'on choisit les coefficients $c_{jk}$ entre $-1/2$ et $1/2$, et tels que $P_n-C_n$ soit à coefficients dans $\Z$. Pour tout $x\in E$, on a :

 \medskip

 $|C_n(x)| \leqslant \frac{1}{2}\sum_{j,k} |x|^j2\lambda^k \leqslant A\sum_{k=0}^{k=n-\l-1}\lambda^k \leqslant A\frac{\lambda^{^{n-\l}}-1}{\lambda-1}   <   A\frac{\lambda^{^{n-\l}}}{\lambda-1}.$

\medskip
\n D'où $|C_n(x)| < 2\lambda^n$, d'après (2.7.8). 

\medskip

\n {\it Démonstration de la partie {\rm(ii)} du th.2.7.1.}

\smallskip
 
   Soient $E_j$ les composantes connexes de $E$, et soit $r_j$ le nombre de racines de $P$ dans $E_j$;
   d'après (2.4.7), le nombre de racines de $Q$ dans l'intérieur de $E_j$ est $r_j-1$. Comme on a supposé que $Q=1$
   (ce qui n'avait pas servi dans la démonstration de (i)), on a donc $r_j=1$, et comme  $\sum r_j = r$, le nombre des $j$ est égal à $r$.
   
    Fixons $j$. On a construit dans la démonstration du 
  th.2.4.5 (voir notamment la formule (2.5.2)) une fonction continue $\theta : E_j\to \R$
  qui donne un homéomorphisme $\theta_j :E_j\to[c_0\pi,c_1\pi]$
  avec $c_1-c_0 = ±r_j$, d'où ici  $c_1-c_0 = ±1$. Choisissons le cas du signe + (l'autre cas est analogue).
  Comme on peut faire une translation sur $c_0$ par un entier, on peut supposer $c_0=0, c_1=1$
  de sorte que $\theta_j$ permet d'identifier $E_j$ à $[0,\pi]$. D'après (2.5.2), cette identification
  transforme le polynôme $P(x)$, pour $x\in E_j$, en $2\lambda
  \cos t$, où $t=\theta(x)$ et elle transforme
  le polynôme $P_n$ du §2.7 en $2\lambda^n\cos nt$. En particulier, les maxima et minima de $|P_n|$ dans $E_j$ correspondent
  à $t = 0, \frac{\pi}{n}$$,\frac{2\pi}{n},...,\pi$, et ses racines correspondent à $t = \frac{\pi}{2n}$$, \frac{3\pi}{2n}, ..., \frac{(2n-1)\pi}{2n}$. 
  Les sous-intervalles de $E_j$ relatifs à $P_n$ correspondent à $[0,\frac{\pi}{n}], [\frac{\pi}{n}, \frac{2\pi}{n}],...,[\frac{(n-1)\pi}{n},\pi]$.
  D'après le lemme 2.5.5, appliqué à $P_n$, le polynôme $P'_n$ choisi ci-dessus a une racine (et, du coup, une seule) dans chacun de ces sous-intervalles. 
  
  \medskip
  
  \n {\bf Lemme 2.7.10.}  {\it Soit $\mu_j$ la restriction à $E_j$ de la mesure d'équilibre $\mu_E$ de~$E$.}
  
  \smallskip
  (a) {\it L'homéomorphisme $\theta_j : E_j \to [0,\pi]$ défini ci-dessus transforme $\mu_j$ en la mesure $\frac{1}{r\pi}dt$
  de} $[0,\pi]$.
  
  \smallskip
  (b){\it Quand $n \to \infty$, la mesure définie par les racines de $P_n$ sur $E_j$ tend vers~$r\mu_j$, et il en est de même
  de la mesure définie par les racines de $P'_n$ sur $E_j$.}
 
 \smallskip
 
 (c)  \ $\mu_E = \lim \delta_{P_n} = \lim \delta_{P'_n}$.
 
 \smallskip
 
 [Il est clair que (c) entraîne la partie (ii) du th.2.7.1.]
 
 \smallskip
 
 \n {\it Démonstration de} (a).
 
  Le composé $\varphi$ de $\theta_j : E_j \to [0,\pi]$ et de $M\!\cos : [0,\pi] \to [M,-M]$
 n'est autre que la restriction de $P$ à $E_j$. En appliquant à $P: E \to [-M,M]$ la formule (A.6.9) du §A.6, on voit
 que l'image de $r\mu_j$ par $\varphi$ est égale à la mesure d'équilibre de $[M,-M]$, laquelle n'est autre que 
 l'image par $M\!\cos$  
 de $\frac{1}{\pi}$$dt$. Cela entraîne~(a).
 
 \smallskip
 
 \n{\small {\it Remarque.} L'assertion (a) peut se reformuler en disant que la mesure d'équilibre de $E$ est la mesure
 $\frac{1}{\pi}$$|\eta|$ où $\eta$ est la forme différentielle définie au §2.3. C'est là un résultat général, valable même pour des $E$ qui ne sont pas de type Pell-Abel; il est souvent cité (par exemple dans [Wi 69], §14, p.226); on en trouvera
 une démonstration dans l'Appendice B, §B.7.\par}

 \smallskip
 \n {\it Démonstration de} (b). Vu (i), cela revient à voir que, quand $n \to \infty$, la mesure définie sur $[0,1]$ par les mesures de Dirac en les points
 $\frac{1}{2n}$$, \frac{3}{2n}, ..., \frac{2n-1}{2n}$ tend vers la mesure de Lebesgue de $[0,1]$, ce qui est clair (c'est le principe du calcul
 des intégrales de fonctions continues par les sommes de Riemann). Même chose pour les racines des $P'_n$, puisqu'il y en a une dans chaque sous-intervalle.  
  
   \smallskip
 \n {\it Démonstration de} (c). Cela résulte de (b) et du fait que $\mu_E = \sum \mu_j$.

 \bigskip

%
%
\subsection{Comment ramener le cas général au cas de Pell-Abel}

  Fixons $g \geqslant 0$. Soit $U$ l'ouvert de $\R^{2g+2}$ formé des points $(a_0,b_0,a_1,\ldots, b_g)$
  tels que $a_0 < b_0 < a_1 <  \cdots < b_g$. Si $u = (a_0,\ldots,b_g)$ appartient à~$U$, notons~$E_u$
  le compact correspondant:
  
  \smallskip
  
  \hspace{20mm} $E_u = [a_0,b_0] \cup \cdots \cup [a_g,b_g].$
  
  \smallskip
  
  Disons que $u$ {\it est de type PA} si l'équation de Pell-Abel correspondante a une solution de degré $>0$.
  Soit $U_{PA}$ l'ensemble des $u$ de ce type.
  
  \smallskip
  
  \n {\bf Théorème 2.8.1.} {\it $U_{PA}$ est dense dans $U$.}
  
  \smallskip
  
  Cet énoncé permet d'achever la démonstration du th.~1.6.1. En effet, soit $u$
  un élément de $U$ tel que $\cp(E_u) > 1$.
   D'après la prop.~A.7.1, tous les $v \in U$ suffisamment voisins de $u$ sont tels que $\cp(E_v) > 1$. D'après le th.~2.8.1, on peut choisir $v\in U_{PA}$  avec $E_v \subset E_u$ et $\cp(E_v) > 1$. Grâce au cor.~2.6.3, on peut aussi supposer que l'équation de Pell-Abel $P^2-DQ^2=M^2$ pour $E_v$ 
  a une solution  avec $P\in \Q[x]$ et $M \in \Q$; le cor.~2.7.2 montre alors que $\Irr_{E_v}$ est infini,
 donc aussi $\Irr_{E_u}$, ce qui démontre le th.~1.6.1.
 
\medskip

\n {\it Démonstration du th.~2.8.1} (d'après [La 16]).

\smallskip

  Si $u\in U$, soit $J_u$ la jacobienne de la courbe hyperelliptique $C_u$ associée à $E_u$. La composante neutre $J_u(\R)^0$
  de $J_u(\R)$ est un tore de dimension $g$, que l'on peut identifier à $\R^g/\Z^g$, après un choix convenable (précisé dans [La 16]) de bases pour les formes de première espèce, et pour les cycles réels. L'image de $\infty_- - \infty_+$
  dans $J_u(\R)$ appartient à~$J_u(\R)^0$; d'où un élément $\vartheta(u) \in \R^g/\Z^g$. Cela définit une application  
  analytique réelle (cf. [La 16]) $\vartheta : U \to \R^g/\Z^g$; comme  $U$  est simplement connexe, on peut la relever en 
  une application continue  $\theta : U \to \R^g$. On a :
  
  \smallskip

  (2.8.2) {\it $u \in U_{PA} \Longleftrightarrow \theta(u) \in \Q^g$.}

\smallskip
  
\n En effet, $\theta(u) \in \Q^g \Longleftrightarrow \vartheta(u)$ est d'ordre fini dans $\R^g/\Z^g$.

\smallskip

Ainsi, le th.~2.8.1 équivaut à dire que {\it l'image réciproque de $\Q^g$ par $\theta : U \to \R^g$ est dense dans~$U$}.
On va voir que cela provient simplement du fait que  $\Q^g$ est dense dans~$\R^g$. Il faut d'abord «rappeler» quelques faits élémentaires :

\medskip

\n {\it Intermède topologique.}

\smallskip
\n {\bf Proposition 2.8.3.} {\it Soit $f : X \to Y$ une application entre espaces topologiques. Soit  $Y'$ une partie dense
de $Y$. Si $f$ est une application ouverte, $f^{-1}(Y')$ est dense dans $X$.}

\n [Rappelons que $f$ est dite {\it ouverte} si l'image par $f$ de tout ouvert de $X$ est un ouvert de $Y$; cela n'entraîne pas que  $f$  soit continue.]

\smallskip
\n {\it Démonstration.} Soit $V$ un ouvert non vide de~$X$. Alors $f(V)$ est un ouvert non vide de~$Y$; il rencontre
donc~$Y'$, d'où $V \cap f^{-1}(Y') \neq \varnothing$, ce qui montre que $f^{-1}(Y')$ est dense dans $X$.
   \smallskip
   
   Nous allons appliquer ceci aux variétés analytiques réelles :
   
    \smallskip
    
    \n {\bf Proposition 2.8.4.} {\it Soit $f: X \to Y$ un morphisme de variétés analytiques réelles de dimension finie.
    Supposons que $X$ soit connexe, et qu'il existe $x \in X$ en lequel $f$ est une submersion. Soit $Y'$ une 
    partie dense de $Y$. Alors $f^{-1}(Y')$ est dense dans $X$.}

\n [Rappelons que $f$ est une submersion en $x$ signifie que l'application tangente à $f$ en $x$ est surjective.]

\smallskip
\n {\it Démonstration.}  Soit $F$ l'ensemble des points de $X$ en lesquels $f$ n'est pas une submersion;
c'est un sous-ensemble analytique fermé de $X$ : il est défini localement par l'annulation d'un nombre fini 
d'équations analytiques. Soit $\overset {\;\circ} F$ l'intérieur de $F$. Le classique «principe du prolongement analytique»
dit que $\overset {\;\circ} F$ est fermé. Comme il est ouvert, et que $X$ est connexe, c'est, soit $\varnothing$, soit $X$. Or ce n'est pas $X$, puisqu'il ne contient pas~$x$. C'est donc $\varnothing$, autrement dit $X \sm F$ est dense dans $X$. Comme la restriction de~$f$ à~$X \sm F$ est une submersion, c'est une application ouverte. La prop.~2.8.3 entraîne que $(X \sm F) \cap f^{-1}(Y')$ est dense dans $X\sm F$, donc aussi dans~$X$.

\medskip\medbreak

\n {\it Fin de la démonstration du th.~2.8.1}
 \par\nobreak
 \smallskip\nobreak
  Pour démontrer le th.~2.8.1, il suffirait, d'après la prop.~2.8.4, d'exhiber un élément~$u$ de~$U$, en lequel l'application
$\theta : U \to \R^g$ est une submersion. Malheureusement, ce n'est pas facile pour $g > 1$. Dans [La 16], Lawrence procède autrement : il définit une variété $U'$ qui contient $U$, ainsi que certains points «dégénérés»; l'application 
$\theta$ se prolonge à $U'$, et il montre qu'elle est une submersion en certains des points dégénérés. Je renvoie à
[La 16] pour plus de détails. {\it Grosso modo}, les points dégénérés qu'il utilise correspondent à 
des suites $a_0 < b_0=a_1 < b_1 = a_2 < \cdots < b_{g-1} = a_g < b_g$. La courbe hyperelliptique de genre $g$ est remplacée par une courbe de genre $0$ ayant $g$ points doubles; sa jacobienne généralisée est un groupe de type multiplicatif de dimension $g$ dont le groupe des points réels est compact.

\medskip

\n {\it Remarque.} On trouvera d'autres démonstrations du th.2.8.1 dans [Ro 64], [Bo 99] et [ACZ 18].

%
%
\subsection{Fin de la démonstration du théorème 1.6.2}

  \smallskip
  
    Conservons les notations du § 2.8 ci-dessus. Soit $E = E_u$, avec $u\in U$, une réunion finie d'intervalles fermés telle que $\cp(E)>1$.
   D'après la prop.~A.8.1, tous les $v \in U$ suffisamment voisins de $u$ sont tels que $\cp(E_v) > 1$. D'après le th.~2.8.1, on peut choisir $v\in U_{PA}$  avec $E_v \subset E_u$ et $\cp(E_v) > 1$. Grâce au cor.~2.6.3, on peut aussi supposer que l'équation de Pell-Abel $P^2-DQ^2=M^2$ pour $E_v$ 
  a une solution  avec $P\in \Q[x], Q=1,M \in \Q$ et $M > 1$. D'après le th.2.7.1, on peut choisir des polynômes $P_{n,v}$ unitaires à coefficients entiers, dont toutes les racines sont dans $E_v$, et qui sont tels que $\delta_{P_{n,v}}$ converge vers $\mu_{E_v}$.
  Choisissons maintenant une suite d'éléments $v_j$ de $U$ ayant les propriétés ci-dessus, et tels que les $E_{v_j}$ forment une suite croissante, d'adhérence égale à $E$. D'après (A.5.5), on a $\mu_E = \lim \mu_{E_{v_j}}$; le procédé diagonal montre qu'il existe des $n_j$ tels que
  $\mu_E = \lim \delta_{P_{n_j,v_j}}$. Cela achève la démonstration du th.1.6.2, et donc aussi celle du th.1.6.1.

%
%
%
\section*{Appendice A. Fascicule de résultats sur les capacités}
\stepcounter{section}\def\thesubsection{A.\arabic{subsection}}

  Soit $K$ une partie compacte de $\C$. La {\it capacité} $\cp(K)$ de $K$ (parfois appelée {\it capacité logarithmique}, ou bien {\it diamètre transfini}) est un nombre réel $\geqslant 0$, défini de l'une des trois façons équivalentes A.1, A.2, A.3  données ci-dessous. Cette notion a été introduite en 1923-1924 par Fekete ([Fe 23]) et  Szeg\H o ([Sz 24]), sans doute inspirés par des résultats antérieurs de Stieltjes ([St~85]) et de Schur ([Sc 18]). On en trouvera une étude détaillée dans Tsuji ([Ts 59], chap.~III) et Ransford ([Ra 95], chap.~3-4-5).

\addtocounter{footnote}{-1}
\subsection{La capacité définie au moyen de discriminants\protect\footnotemark}
\footnotetext{Dans les énoncés ci-dessous, on suppose  $K$  non vide;
 lorsque $K$ est vide, et plus généralement quand $K$ est fini, on a  $\cp(K)=0$.}

 Si $n > 1$, posons:
 
 \smallskip
 
 (A.1.1) \  $d_n(K) = \sup_{x_1,\ldots,x_n \in K} \prod_{i\neq j} |x_i-x_j|^{1/n(n-1)}.$
 
 \smallskip
 
\n{\small [Lorsque $n=2$, $d_n(K)$ est le {\it diamètre} de $K$, au sens habituel.]\par}

 \smallskip

 On a \ $d_2(K) \geqslant d_3(K) \geqslant \cdots$. Lorsque $K$ est fini, on convient que $d_n(K)=0$
 si $n > |K|$.
 
 \smallskip
 
 \n La capacité de $K$ est définie par :
 
 \smallskip
  (A.1.2) $ \ \cp(K) = \inf_n d_n(K) = \lim_{n \to \infty} d_n(K)$.
  
  \smallskip
  
\n   Elle est souvent notée $d_\infty(K)$. 
 
 
\subsection{La capacité, à la Chebyshev}
 
   Si $n > 0$, soit 
$c_n(K) = \inf_P\mathopen\|P\mathclose\|_K^{1/n}$, 
où $P$ parcourt l'ensemble des polynômes unitaires de degré  $n$
   à coefficients dans $\C$, et 
$\mathopen\|P\mathclose\|_K$ 
est la borne supérieure de $|P|$ sur $K$. Pour chaque $n \leqslant \card(K)$, il existe un unique $P$ (cf. [Ts 59], th.~III.23) tel que 
$\mathopen\|P\mathclose\|_K = c_n(K)^n$: 
c'est le {\it $n$-ième polynôme de Chebyshev} de $K$.
      
    \smallskip
   On a :      
 \smallskip

   (A.2.1) \  \quad $\cp(K) = \inf_n c_n(K)^{1/n} = \lim_{n \to \infty} c_n(K)^{1/n}.$

   \smallskip
   
   \n {\it Exemple.} Prenons pour $K$ le segment $[-2,2]$; il est bien connu que le $n$-ième polynôme de Chebyshev
   de $K$ est le polynôme $T_n$ caractérisé par
   
    \smallskip

   (A.2.2) \ $T_n(t+t^{-1}) = t^n + t^{-n}$,
 
  \smallskip
  
 \n  ou, ce qui revient au même :
   
    \smallskip

    (A.2.3) \  $T_n(2 \cosh x) = 2 \cosh nx$  \ et \ $T_n(2 \cos \theta) = 2 \cos n \theta$. 
    
     \smallskip

    On a $\parallel\!T_n\!\parallel_K = 2$; d'après (A.2.1), cela entraîne $\cp(K) = 1$. Ce résultat peut aussi se
    déduire de (A.1.2), et de la détermination des  $d_n(K)$ due à Stieltjes et Schur, cf. [St~85] et [Sc 18], §1, Satz I.
    
    La mesure d'équilibre de $K$ (au sens de A.4 ci-dessous) est :
    
    \smallskip
    (A.2.4)  \ $\mu_K = \frac{1}{\pi}\frac{dx}{\sqrt{4-x^2}}$; 
    
    \smallskip
    \n cela résulte de la prop.A.7.1, et du fait que la mesure d'équilibre d'un cercle est
    l'unique mesure de masse 1 invariante par rotation.
   
   \medskip
  Par homothétie, on déduit de  $\cp([-2,2]) = 1$
    que {\it la capacité d'un intervalle de longueur  $\l$  est égale à} $\l/4$.
   
   \medskip

%
%
%
\subsection{Une variante de A.1, en termes de mesures}

       Soit $\mu$ une mesure positive sur  $K$. Posons :
       
        \smallskip
       
     (A.3.1) \ $I(\mu) = \iint_{K \times K} \log|x-y|\mu(x)\mu(y).$
       
        \smallskip
       
   \n   C'est, soit $-  \infty$, soit un nombre réel. Soit $v(K) = \sup_\mu I(\mu)$, où $\mu $ parcourt l'ensemble des mesures positives de masse 1 à support dans~$K$. On a :
   
    \smallskip
     
  (A.3.2) \  $\cp(K) = e^{v(K)}$ ,
     
      \smallskip
     
  \n et en particulier  $\cp(K) = 0$ si et seulement si  $I(\mu) = -\infty$ pour tout~$\mu$ à support dans~$K$.
  
   \smallskip

   {\small Le fait que les définitions A.1, A.2 et A.3 de $\cp(K)$ sont équivalentes est dû à Fekete et Szeg\H o, cf. [Fe 23], [Sz 24], [FS 55], ainsi que [Ra 95], th.~5.5.2 et th.~5.5.4. La terminologie « diamètre transfini » provient de A.1 et celle de « capacité logarithmique » de A.3. Noter que, du point de vue de la théorie des capacités de Choquet, c'est $v(K) = \log \cp(K)$, et non $\cp(K)$, qui mériterait le nom de « capacité », cf. [Ch 58].\par}

\medskip

\n  A.3.3. {\it Extension de la notion de capacité aux ensembles non compacts.}  

Soit $Y$ une partie bornée de~$ \C$. On appelle {\it capacité intérieure} (ou simplement « capacité ») de $Y$, la borne supérieure des $\cp(K)$ lorsque $K$ parcourt les sous-espaces compacts de $Y$ ; on la note encore $\cp(Y)$. En particulier, $Y$ est dit  {\it de capacité $0$} («ensemble polaire» : « polar set » dans [Ra 95]) si $\cp(K) =0$ pour toute partie compacte~$K$ de~$Y$.

%
\subsection{Mesure d'équilibre}

  Lorsque $K$ est un compact de capacité $>0$, il existe une unique mesure positive $\mu$ de masse 1 telle que $I(\mu) = \log \cp(K)$, cf. [Ra 95], th.~3.7.6. On l'appelle la {\it mesure d'équilibre} de $K$, et on la note $\mu_K$. C'est une mesure diffuse (cela résulte du cor.~1.4.2). Son support
 Supp$(\mu_K)$ n'est pas toujours égal à $K$; ainsi, lorsque  $K$  est un disque, le support de $\mu_K$ est le cercle qui borde ce disque. Lorsque $K$ est contenu dans $\R$, $K$ et $\Supp(\mu_K)$ ne diffèrent que par un ensemble de capacité $0$. De façon plus précise:
 
 \smallskip
 
 \n {\bf Proposition A.4.1.} {\it Soit $K$ une partie compacte de $\R$ de capacité $>0$. Soit $K'$ une partie fermée de $K$. Les propriétés suivantes sont équivalentes} :
 
  \smallskip
  (A.4.2) $ \cp(K') = \cp(K)$.
  
   \smallskip
  (A.4.3) $ K' \supset \Supp(\mu_K).$
  
   \smallskip
  (A.4.4) $K \sm K'$ {\it est de capacité}~$0$ (au sens de A.3.3).
 
  \medskip
 
  \n {\bf Corollaire A.4.5.} {\it $\Supp(\mu_K)$ est la plus petite partie fermée de $K$ ayant même capacité que $K$.}

  \medskip
 
  \n {\bf Corollaire A.4.6.} {\it L'ensemble 
$K \sm \Supp(\mu_K)$ 
est de capacité~$0$.}

 \medskip\medbreak
 
 \n {\it Démonstration de la prop.~A.4.1.}
 \par\nobreak
 \smallskip\nobreak
\n (A.4.2) $\Rightarrow$ (A.4.3). Si $\cp (K') = \cp (K)$, on a $I(\mu_{K'}) = I(\mu_K)$, d'où $\mu_{K'} = \mu_K$ en vertu de l'unicité de la mesure d'équilibre de $K$; cela entraîne $\Supp(\mu_{K'}) = \Supp(\mu_K)$, d'où $ K' \supset \Supp(\mu_K).$

 \smallskip
\n  (A.4.3)  $\Rightarrow$ (A.4.4). D'après [Ts 59], th.~III.31, l'ensemble  
$K \sm \Supp(\mu_K)$ 
est de capacité~$0$.
Il en est a fortiori de même de 
$K \sm K'$ 
si $K'$ contient $\Supp(\mu_K)$.

\smallskip

 \n (A.4.4) $\Rightarrow$ (A.4.2). Cela résulte de A.5.4 ci-dessous.

   \medskip
  
  \n A.4.7. Disons que $K$ est {\it réduit} si $K = \Supp(\mu_K)$. La prop.~A.4.1 entraîne:
  
  \medskip
  
  \n {\bf Proposition A.4.8.} {\it Soit $K$ une partie compacte de $\R$ de capacité~$>0$. Les propriétés suivantes sont équivalentes} :
  
  (i) \  $K$ {\it est réduit.}
  
  (ii) {\it Aucune partie fermée de $K$, distincte de $K$, n'a la même capacité que}~$K$.
  
  (iii) {\it Aucune partie ouverte non vide de $K$ n'est de capacité $0$.}

  \medskip
  
  \n {\bf Corollaire A.4.9.} {\it Soit $\mu$ une mesure positive à support compact sur $\R$, telle que $I(\mu) > - \infty$.
  Alors $\Supp(\mu)$ est réduit.}
    
  \smallskip
  
  \n {\it Démonstration.} Soit $K = \Supp(\mu)$. Si $K$ n'était pas réduit, d'après la prop.~A.4.8, il existerait 
  un ouvert non vide $U$ de $K$  de capacité 0. D'après le th.~III.7 de [Ts 59], on aurait $\mu(U) = 0$,
  ce qui contredirait le fait que $U$ est contenu dans $\Supp(\mu)$.
    
%
%
\subsection{Quelques propriétés de la capacité}
  
  (A.5.1) ({\it Linéarité}) \ $\cp(\lambda K) = \mathopen|\lambda\mathclose| \cp(K)$ \ pour tout $\lambda \in \C$.

 \smallskip
         
(A.5.2) ({\it Continuité pour les suites décroissantes}) Soit $K_n$ une suite décroissante de compacts de $\C$. On a $\cp(\cap K_n) = \inf_n  \cp(K_n)$, cf. [Ra 95], th.~5.1.3 (a).

 \smallskip
 
 (A.5.3) Soient $K_1, K_2 \subset \C$ deux compacts, et soit $d$ leur distance. On a :
 
  \smallskip
 \hspace{5mm}  $\cp(K_1 \cup K_2) \ \geqslant \ \cp(K_1)^{1/4}\cp(K_2)^{1/4}d^{1/2}$.
  
   \smallskip 
   
   Cela se démontre en appliquant (A.1.1) avec $n$ pair, en choisissant de façon optimale  $n/2$ points dans $K_1$ et $n/2$ points dans $K_2$.
   
   \smallskip
      
      (A.5.4) Soient $B_1$ et $B_2$ deux parties boréliennes bornées de $\C$. Si $B_2$ est de capacité~$0$, on a $\cp(B_1 \cup B_2)=\cp(B_1),$ cf. [Ts 59], th.~III.18.
      
          \smallskip
      
      (A.5.5) ({\it Continuité pour les suites croissantes}) Soit $K_n$ une suite croissante de compacts et soit $K$ un compact contenant $\cup_n K_n$. Supposons que $\cp(K \sm \cup K_n) = 0.$ Alors $\cp(K) = \sup \cp(K_n).$

      \smallskip
      \n {\it Démonstration}. 
      
      D'après [Ra 95], th.~5.1.3 (b), on a $\cp(\cup_n K_n) = \sup \cp(K_n)$. On applique (A.5.4)
      à $B=K, B_1= \cup_n K_n$ et $ B_2 = K \sm \cup_n K_n.$

      \smallskip
         (A.5.6) Soit $K$ un compact contenu dans $\R$, et soit mes$(K)$ sa mesure de Lebesgue.  On a $\cp(K)\geqslant \  $mes$(K)/4$, cf.  [Ra 95], th.~5.3.2.
         
         \smallskip

En particulier,  $\cp(K)=0$ entraîne mes$(K)=0$. La réciproque est fausse: l'ensemble triadique de Cantor dans $[0,1]$ est de mesure nulle, mais sa capacité est au moins $1/9$, cf. §1.6.
         
 \smallskip
         
(A.5.7) Soit $f \in \C[X]$ un polynôme unitaire de degré $d \geqslant 1$.

\n On a   $\cp(f^{-1}K) = \cp(K)^{1/d}$ pour tout compact $K$ de $\C$, cf. [Ra 95], th.~5.2.5.

    \smallskip
   
    (A.5.8) ({\it Capacité de la réunion de deux intervalles de même longueur})
    
     \smallskip

    Soient  $a,b \in \R$ avec $0<a <b$, et soit  $E = [-b,-a] \ \cup \ [a,b]$. 
  On a : 

\smallskip
 \hspace{5mm} $\cp(E) = \frac{1}{2}\sqrt {b^2-a^2}.$

 Cela résulte
    de (A.5.7) appliqué à $f = X^2$ et $ K = [a^2,b^2]$, de sorte que $f^{-1}K =E$.

%
%
\subsection {\bf Capacité des images réciproques}

   \smallskip
     Soit $f \in \C[X]$ un polyn\^{o}me unitaire de degré $d \geqslant 1$. Soit $K$ un compact de $\C$
     et soit $L = f^{-1}(K)$. Comme on l'a dit ci-dessus (cf. (A.5.7)), on a :
     
     \smallskip
     
     (A.6.1) \ $\cp(L) = \cp(K)^{1/d}$.
     
     \smallskip

     Supposons $\cp(K) > 0$, de sorte que $\cp(L) > 0$; les mesures d'équilibre $\mu_K$ et $\mu_L$ 
     sont donc bien d\'{e}finies. Voici comment on peut passer de l'une à l'autre :
     
     \smallskip
     
     \n {\bf Proposition A.6.2}. (i) $\mu_K $ {\it est l'image de $\mu_L$ par $f: L \to K$.}
     
     (ii) $\mu_L $ {\it est l'image réciproque de $\mu_K$ par $f:L \to K$, au sens défini ci-dessous.}
   
   \medskip

     \n {\bf Définition de l'image réciproque $\nu^*$ d'une mesure $\nu$ sur $K$.}
     
     \smallskip
   Il faut d'abord définir une opération d'{\it image directe} pour les fonctions. Soit $\varphi \in C(L)$.
   Soit $x \in K$; l'image réciproque de $x$ par $f$ est l'ensemble $x_1,...,x_d$ des racines du polyn\^{o}me
   $f(z)-x$; on convient de répéter  chaque racine d'après sa multiplicité. Cela donne un sens
   à l'expression  $\varphi_*(x) = \frac{1}{d}$$(\varphi(x_1) + \cdots + \varphi(x_d))$. On constate facilement 
   que $\varphi_*$ est une fonction continue sur $K$ (c'est clair en dehors des valeurs critiques de $f$;
   il faut un argument local pour les valeurs critiques). On obtient ainsi une application linéaire
   continue $C(L) \to C(K)$. On en déduit par dualité une application linéaire sur les mesures : à une mesure $\nu$
   sur $K$, on associe la mesure $\nu^*$ sur $L$ telle que :
  
  \smallskip
   
     (A.6.3)  $ \nu^*(\varphi) = \nu(\varphi_*)$ \ \ pour tout  $\varphi \in C(L)$. 
     
    \smallskip  
   
   Si $\nu$ est positive de masse 1, il en est de même de $\nu^*$; de plus, l'image de $\nu^*$
   par $f$ est égale à $\nu$. On peut voir $\nu^*$ comme le {\it relèvement canonique} de $\nu$.
   
   \medskip
   \n {\small [Variante: si $x\in \C$, les $d$ racines $x_1,...,x_d$ de l'équation $f(X)=x$ définissent 
  une mesure $\delta(x) =  \frac{1}{d}$$\sum \delta_{x_i}$ qui dépend contin\^{u}ment de $x$,  et l'on a  $\nu^* = \int_K \delta(x)\nu(x)$, l'intégrale étant prise dans l'espace des 
  mesures, comme dans [INT], chap.V, §4.]\par}
  
  \medskip
  
  \n {\it Remarque.} Plus généralement, on peut prendre pour fonction $\varphi$ n'importe quelle fonction continue sur $L$ à valeurs dans $\R \cup \{-\infty\}$; la fonction $\varphi_*$ correspondante jouit de la même propriété, et la formule (A.6.3) est encore valable : c'est immédiat en tronquant et en passant à
  la limite. Nous en aurons besoin ci-dessous.
   
   \medskip
   \n {\bf Une première formule}
   
   \smallskip
   
     Soient $x,y \in \C$, soient $x_1,...,x_d$ les $d$ racines de l'équation $f(X)=x$, et soient $y_1,...,y_d$ celles de $f(X)=y$.
     Alors :
     
     \smallskip
     
     (A.6.4) \  $ \sum_{i=1}^d \log |x_i - y_j| = \log|x - y|$  pour tout $j = 1,...,d$.
     
     \medskip
     
     \n Cela entra\^{i}ne :
     
     \smallskip
     
     (A.6.5) \ $ \frac{1}{d}$ \!$\sum_{i,j} \log |x_i - y_j| = \log|x - y|$.
     
     \medskip
   
   \n {\it Démonstration.} Les $x_i$ sont les racines du polynôme unitaire $f(X)-x=0$. On a donc :

    \smallskip
    
    (A.6.6) \  $f(X)-x = (X-x_1)\cdots (X-x_d)$.
    
    \smallskip
    
    \n En remplaçant $X$ par $y_j$, et en prenant les valeurs absolues des deux membres, on obtient  $|y-x| = \prod |y_j-x_i|$, d'où (A.6.4).

   \bigskip
  
  \n {\bf Une seconde formule}
  
  \smallskip
  
   Soit $\nu$ une mesure positive sur $K$, et soit $\nu^*$ la mesure correspondante sur~$L$. On a :
   
   \smallskip
   
   (A.6.7)  \ $  I(\nu^*) =  \frac{1}{d}$$I(\nu). $
  
  \smallskip
  
 \n  [Rappelons que :
  
   \smallskip

   $I(\nu) = \iint_{K \times K} \log |x-y| \nu(x)\nu(y))$ \ et \ $I(\nu^*)=\iint_{L \times L} \log |u-v| \nu^*(u)\nu^*(v)$.] 
    
   \medskip
   
   \n {\it Démonstration de} (A.6.7). 
   
    \smallskip
    
 \n   Posons, pour simplifier l'écriture:
    
     \smallskip
(A.6.8) \ $ C(x,y) = \log|x-y|.$

 \smallskip
\n Soit $h$ la fonction définie par :

 \smallskip
  $h(u) = \ \int_L C(u,v) \nu^*(v)$.
  
   \smallskip
 \n Par définition de $\nu^*$, on a :
  
   \smallskip
  $h(u) = \int_K C_*(u,y) \nu(y)$, \ où \   $C_*(u,y) = \frac{1}{d}$ \!$\sum_j C(u,y_j)$,
  
   \smallskip
\n  les $y_j$ étant les racines de  $f(X)=y$.
  
   \smallskip
 \n D'après le théorème de Lebesgue-Fubini ([INT], chap.V, §8, prop.5), on a
    
     \smallskip
    $I(\nu^*)= \int_L h(u)\nu^*(u) =  \frac{1}{d}$ \!$\iint_{L \times K} \sum_j C(u,y_j) \ \nu^*(u)\nu(y)$.
    
     \smallskip
  \n  Le même argument que pour $h$ montre que ceci s'écrit aussi, en termes des solutions $x_i$
    de $f(X)=x$, comme :
    
     \smallskip
    $I(\nu^*) =  \frac{1}{d^2}$$ \iint_{K \times K} \sum_{i,j} C(x_i,y_j) \ \nu(x)\nu(y)$.
    
     \smallskip
 \n   D'après (A.6.5), cela donne:
    
     \smallskip
    $I(\nu^*) =  \frac{1}{d}$$ \iint_{K \times K} C(x,y) \ \nu(x)\nu(y) =  \frac{1}{d}$$I(\nu)$.
        
    \medskip
    
    \n {\bf Démonstration de la prop.A.6.2}
    
    \smallskip 
  \n  On applique (A.6.7) à $\nu = \mu_K$. On en déduit :
    
    \smallskip
    \hspace{15mm} $I(\mu_K^*) = \frac{1}{d}$ \!$I(\mu_K) = \frac{1}{d}$ \!$ \log \cp(K) = \log \cp(L)$,
    
    \smallskip
 \n    ce qui prouve que $\mu_K^ *$ est la mesure d'équilibre de $L$; il est clair que son image par $f$ est $\mu_K$.  
    
    \medskip
    
    \n {\bf Un cas particulier.}
    
    \medskip
      C'est celui où $L$ est réunion de $d$ parties compactes disjointes $L_i \ (i = 1,...,d)$ telles que les projections $: L_i \to K$
      définies par $f$ soient des homéomorphismes. Soit $\nu$ une mesure sur $K$; notons $\nu_i$ la mesure sur $L_i$ déduite de 
      $\nu$ par l'homéomorphisme $K \to L_i$ inverse de $f$; on peut considérer $\nu_i$ comme une mesure sur~$L$. On a:
      
      \smallskip 
      
      (A.6.9) \ $ \nu^* =    \frac{1}{d}$$(\nu_1 + \cdots + \nu_d).$
      
      \smallskip
      
    \n   C'est immédiat.
      
      \medskip
      
      \n Une autre façon d'exprimer (A.6.9) consiste à identifier $L$ au produit $K \times \{1,...,d\}$, et à écrire $\nu^*$ comme un produit tensoriel :
      
      \smallskip
      (A.6.10) \ $\nu^* = \nu \otimes  \frac{1}{d}$$(\delta_1 + \cdots + \delta_d)$.

%
%
\subsection{Capacité des sous-ensembles d'un cercle}

L'énoncé suivant est analogue à ceux de A.6: 

    \smallskip
    
    \n {\bf Proposition A.7.1.} {\it Soit $C$ un cercle de centre $0$ et de rayon $r$. Soit $I = [-2r,2r]$, et soit
    $f : C \to I$ l'application  $z \mapsto z + \zbar.$ Soit $K$ un compact de $I$ et soit} $K_C = f^{-1}(K)$.
    
     (i) {\it On a}  $\cp(K_C) = r^{1/2} \cp(K)^{1/2}.$
     
     (ii) {\it Supposons} $\cp(K) \neq 0$. {\it On a $f(\mu_{K_C}) = \mu_K$, où  $\mu_K$ et $\mu_{K_C}$ sont les mesures d'équilibre
     de  $K$ et de $K_C$, cf.~A.4}. 

  \smallskip
  
  \n{\small La notation $f(\mu_{K_C})$ désigne {\it l'image} de la mesure $\mu_{K_C}$ par l'application $f$.\par}
  
  \smallskip
  
  \n {\it Remarque.} L'assertion (i) se trouve déjà dans [Ro 69].
  
  \medskip
  \n {\bf Corollaire A.7.2.} {\it La capacité d'un cercle de rayon  $r$  est égale à $r$.}
    
    \smallskip
    
    Cela résulte de (i) appliqué à $K= I, K_C = C$.
    
    \smallskip

 \n {\it Démonstration de la prop.~A.7.1.}
 
 \smallskip
 
 Quitte à faire une homothétie, on peut supposer $r=1$; alors $C$ est le cercle unité $|z|=1$,
  et $I = [-2,2]$. On peut aussi supposer que $\cp(K) > 0$.
  
    Soit $\nu$ une mesure positive sur $C$, invariante par la conjugaison complexe  $z \mapsto \zbar$. Soit~$\nu'$ l'image de~$\nu$ par~$f$; c'est une mesure sur~$I$, de même masse que~$\nu$. Soient $I(\nu)$ et $I(\nu')$ les intégrales définies
dans (A.3.1), autrement dit :

\smallskip
    
      $I(\nu) = \iint_{C \times C} \log|z_1-z_2|\nu(z_1)\nu(z_2)$
    
   \n   et 
     
     $I(\nu') = \iint_{I \times I} \log|x-y|\nu'(x)\nu'(y).$
    
    \medskip
  
  \n {\bf Lemme A.7.3.} {\it On a}  $I(\nu') = 2 I(\nu)$.
  
   \smallskip

  \n {\it Démonstration du lemme.}
  
  Posons $A(z_1,z_2) = \log|z_1-z_2| + \log |\zbar_1 -z_2|. $ Comme  $\nu$ est invariante par conjugaison, on a
  
  \smallskip

(A.7.4) $\iint_{C \times C} A(z_1,z_2)\nu(z_1)\nu(z_2) = 2 I(\nu).$

 \smallskip
 
\n Soient $x = f(z_1)= z_1+ \zbar_1$ et $y = f(z_2)$. Un calcul simple montre que :
 
  \smallskip
  
  (A.7.5)  $|(z_1 - z_2)(\zbar_1 - z_2)| = |z_1+\zbar_1 - (z_2+\zbar_2)|  = |x-y|$,
  
   \smallskip
 \n d'où:
  
   \smallskip
   
   (A.7.6) $A(z_1,z_2) = \log|x-y|$.
   
   \smallskip
   
   Ainsi, la fonction $A(z_1,z_2)$ est la composée de $f \times f : C \times C \to I \times I$ et de la fonction $\log |x-y|$ sur $I \times I$.
   
    \n Comme $\nu' \otimes \nu'$ est l'image de $\nu \otimes \nu$ par $f \times f$, on en déduit :

   \smallskip
\n (A.7.7)   $\iint_{C \times C} A(z_1,z_2)\nu(z_1)\nu(z_2) = \iint_{I \times I} \log|x-y|\nu'(x)\nu'(y) =~I(\nu')$.
 
 \smallskip
 
 \n Le lemme résulte de (A.7.4) et (A.7.7).
 
 \medskip
 
 \n {\it Fin de la démonstration de la prop.~A.7.1.}
    
    \smallskip
 L'application  $\nu \mapsto \nu'$ donne une bijection entre les mesures positives de masse~1 sur~$K_C$ qui sont invariantes par conjugaison, et les mesures positives de masse 1 sur $K$.

 Soit $\mu_K$ la mesure d'équilibre de $K$; on a $I(\mu_K) = c$, avec
 $c = \log \cp(K)$. Si $\nu$ est la mesure correspondante sur  $K_C$, le lemme A.7.3 montre que $I(\nu) =c/2$,
 d'où $\cp(K_C)  \geqslant  \cp(K)^{1/2}$. Si cette inégalité était stricte, on aurait
  $I(\mu_{K_C}) > c/2$ ; comme la mesure $I(\mu_{K_C})$ est canonique, elle est invariante par conjugaison, et elle correspondrait
 à une mesure $\nu'$ sur $K$ telle que  $I(\nu') > c$, contrairement à la définition de  $c$. On a donc $\cp(K_C) =  \cp(K)^{1/2}$, ce qui démontre (i). On voit en outre que $\nu$ est la mesure d'équilibre de $K_C$, ce qui démontre (ii).
   
%
%
\subsection{Une propriété de continuité pour les réunions finies d'intervalles fermés de $\R$}
  
  Soit $a = \{a_1,\ldots,a_n\}$ une suite strictement croissante de nombres réels, en nombre $n$ pair, et posons :

 \hspace{20mm}  $E_a  =  [a_1,a_2] \cup  [a_3,a_4] \cup  \cdots \cup[a_{n-1},a_n]$. 
   
   \smallskip
   
   \n {\bf Proposition A.8.1.} {\it La capacité de $E_a$ dépend continûment de $a$.}
   
   \smallskip
   
   \n {\it Démonstration.}  Soit $d = \inf  (a_{i+1}-a_i)$. Si $\varepsilon$ est $>0$ et $< d/2$, notons $a'_\varepsilon$ la suite des
   $a_i + (-1)^ i\varepsilon$ et $a''_\varepsilon$ la suite des $a_i - (-1)^ i\varepsilon$. On a:
   
   \smallskip
   
  \quad    $E_{a''_\varepsilon} \ \subset E_a \ \subset E_{a'_\varepsilon}$.
      
      \smallskip
    \n  La continuité de  $\cp(E_a)$ équivaut à dire que $\cp(E_a)$ est la borne inférieure des $\cp(E_{a'_\varepsilon})$
    ainsi que la borne supérieure des $\cp(E_{a''_\varepsilon})$, ce qui résulte respectivement de (A.5.2) et (A.5.5).

%
%
%
\subsection{Calculs effectifs de capacités}
Le calcul de la capacité d'un compact $K$ donné est un problème difficile. On trouvera dans [Ra 95, p.~135] une liste de quelques cas connus. Signalons par exemple celui d'un arc de cercle, de rayon $r$ et d'angle $\alpha$, avec $0 \leqslant \alpha \leqslant  2\pi$ : c'est $r \sin (\alpha/4)$; cela se déduit de la prop.~A.6.1 appliquée à $K=[2r \cos(\alpha/2),2r]$. 
 Voir aussi [RR 07] et [LSN 17] pour des calculs approchés, sur ordinateur.

  \smallskip
  Pour les sous-espaces de $\R$, un cas particulièrement intéressant est celui où  $K$  est réunion disjointe de segments fermés. Le cas de deux segments a été traité par Akhiezer : [Ak 30] et [Ak 32]; le résultat s'exprime en termes de fonctions thêta à la Jacobi. Pour le cas général, voir Bogatyrëv [Bo 99], Peherstorfer-Schiefermayr [PS 99] et Bogatyrëv-Grigoriev [BG 17].
  
\newpage
\centerline{\bf APPENDICE B. POTENTIELS}

\medskip
  
  \centerline{ par Joseph OESTERL\'E}

  \bigskip

\noindent{\bf B.1. Rappels sur les fonctions sous-harmoniques}

\smallskip

Soit $U$ une partie ouverte de $\C$. Une fonction numérique $f:U\to\overline{\R}$ est dite {\it sous-harmonique} si  :\par
$(i)$ elle est à valeurs dans $\R\cup\{-\infty\}$;\par
$(ii)$ elle est semi-continue supérieurement;\par
$(iii)$ pour tout disque fermé $\overline{D}(a,r)$ contenu dans $ U$, $f(a)$ est majoré par la valeur moyenne de $f$ dans $\overline{D}(a,r)$.

\medskip

{\small Sous les hypothèses $(i)$ et $(ii)$, la valeur moyenne de $f$ dans un disque fermé contenu dans $ U$ est un élément bien défini de $\R \cup\{-\infty\}$, cf. 1.1.\par}

\medskip

Lorsque $f$ satisfait les conditions $(i)$ et $(ii)$, la condition $(iii)$ est équivalente à :\\
$(iii')$ Pour tout disque fermé $\overline{D}(a,r)$ contenu dans $ U$, $f(a)$ est majoré par la valeur moyenne de $f$ dans le cercle $C(a,r)$.

De plus, pour qu'elle soit satisfaite, il suffit qu'il existe pour tout $a\in U$ des disques fermés $\overline{ D}(a,r)$ contenus dans $U$ de rayon $r>0$ arbitrairement petit 
pour lesquels $f(a)$ est majoré par la valeur moyenne de $f$ dans $\overline{ D}(a,r)$ (resp. dans $ C(a,r)$).

{\it La propriété pour une fonction numérique dans $ U$ d'être sous-harmonique est donc de nature locale.}

\medskip

Une fonction sous-harmonique $f$ dans $U$ qui n'est égale à $-\infty$ dans aucune composante connexe de $U$ est localement intégrable ([Ra 95], th. 2.5.1). Son laplacien 
 au sens des distributions est une mesure de Radon positive ([Ra 95], th. 3.7.2). Si $g$ est une seconde fonction sous-harmonique  dans $U$ qui n'est égale à $-\infty$ dans aucune composante connexe de $U$ et a même laplacien que $f$, il existe une fonction harmonique $h$ dans $U$ telle que $g=f+h$ ([Ra 95], Lemme 3.7.10).
 
 \bigskip

\noindent{\bf B.2. Potentiels de mesures positives à support compact}

\smallskip

Soit $\mu$ une mesure de Radon positive dans $\C$, à support compact. On appelle {\it potentiel de $\mu$} la fonction $p_\mu:\C\to\R\cup\{-\infty\}$ définie par
$$p_\mu(z)=\int_\C\log|w-z|\mu(w).\leqno{\rm (B.2.1)}$$
Cette fonction est sous-harmonique dans $\C$, harmonique en dehors du support de $\mu$ et l'on a
$$p_\mu(z)=\mu(\C)\log|z|+o(1)\leqno({\rm B}.2.2)$$
lorsque $|z|$ tend vers $+\infty$ ([Ra 95], th. 3.1.2). Le laplacien de $p_\mu$ au sens des distributions est $2\pi\mu$ ([Ra 95], th. 3.7.4).

\medskip

\noindent {\bf Proposition $\bf B.2.3.$} {\it Soit $ K$ une partie compacte de $\C$. L'application $\mu\mapsto p_\mu$ est une bijection de l'ensemble des mesures de Radon positives dans $\C$, de masse $1$ et à support dans $ K$, sur l'ensemble des fonctions $p:\C\to\overline{\R}$ possédant les propriétés suivantes $:$\par
$(i)$ la fonction $p$ est sous-harmonique dans $\C;$\par
$(ii)$ elle est harmonique en dehors de $ K;$\par
$(iii)$ on a $p(z)=\log|z|+ o(1)$ lorsque $|z|$ tend vers $+\infty$.}

Il résulte du début de ce numéro que cette application est bien définie. La relation $\Delta p_\mu=2\pi\mu$ montre qu'elle est injective. Démontrons qu'elle est surjective. Soit donc 
$p:\C\to\overline{\R}$ une fonction possédant les propriétés $(i)$, $(ii)$ et $(iii)$. D'après B.1, $p$~est localement intégrable et son laplacien $\Delta p$ au sens des distributions est une mesure de Radon positive dans $\C$, à support dans $ K$. \'Ecrivons cette mesure $2\pi\mu$. On a alors $\Delta p_\mu=\Delta p$. D'après B.1, il existe 
une fonction harmonique $h$ dans $\C$ telle que $p_\mu=p+h$. On a 
$$h(z)=(\mu(\C)-1)\log|z|+ o(1) \ \leqno ({\rm B}.2.4)$$
et {\it a fortiori} $h(z)=o(|z|)$, lorsque $|z|$ tend vers $+\infty$. Cela implique que $h$ est constante, d'après une variante du théorème de Liouville (cf. [Ra 95], exerc. 1.3.5).  La relation (B.2.4) implique alors d'abord que $\mu(\C)=1$, puis que $h=0$. On a ainsi $p_\mu=p$,  ce qui démontre la surjectivité de l'application considérée.\bigskip

\noindent{\bf B.3. Une propriété des ensembles polaires}

\smallskip
\noindent {\bf Proposition $\bf B.3.1.$} {\it Soit $\mu$ une mesure de Radon positive dans $\C$ à support compact telle que $ I(\mu)>-\infty$. Toute partie borélienne polaire de $\C$ est de mesure nulle pour $\mu$.}\par
Cela résulte de [Ra 95], th. 3.2.3. 

\bigskip

\noindent{\bf B.4. Potentiels des mesures d'équilibre}\smallskip

Soit $K$ une partie compacte de $\C$ non polaire. Notons $\mu_K$ sa mesure d'équilibre  et $p_K$ le potentiel de la mesure $\mu_K$. Rappelons que $v(K)=I(\mu_K)$ est le logarithme de la capacité de $K$ (cf. A.4).

 \medskip

\noindent {\bf Proposition $\bf B.4.1.$} {\it $a)$ On a $p_K\geqslant v(K)$ dans $\C$.\par
$b)$ L'ensemble des points $x\in K$ tels que $p_K(x)>v(K)$ est polaire.}\par
Il s'agit là d'un théorème de Frostman, cf. [Ra 95], th. 3.3.4.

 \medskip

Inversement :

\medskip

\noindent {\bf Proposition $\bf B.4.2.$} {\it Soit $\nu$ une mesure de Radon positive dans $\C$, de masse $1$ et à support dans $K$. Soit $a$ la borne inférieure de $p_\nu$ 
dans $K$. On  a $a\leqslant v(K)$, avec égalité si et seulement si $\nu$ est la mesure d'équilibre de $K$.}\par
On a, d'après le théorème de Lebesgue-Fubini ([INT], chap V, \S\ 8, prop. 5), 
$$I(\nu)=\iint_{K\times K}\log|x-y|\nu(x)\nu(y)=\int_Kp_\nu(y)\nu(y),\leqno({\rm B}.4.3)$$
et par suite $a\leqslant I(\nu) \leqslant v(K)$. Si $a=v(K)$, on a  $I(\nu) = v(K)$, et alors $\nu$ est la mesure d'équilibre de $K$ (cf. A.4).

\bigskip

\noindent{\bf B.5. Caractérisation du potentiel d'une mesure d'équilibre}

\smallskip

\noindent {\bf Proposition $\bf B.5.1.$} {\it Soit $K$ une partie compacte de $\C$, non polaire. Pour qu'une fonction  $p:\C\to\overline{\R}$ soit le potentiel
de la mesure d'équilibre de $K$, il faut et il suffit qu'elle possède les propriétés suivantes :\par
$(i)$ elle est sous-harmonique dans $\C$ et harmonique en dehors de $K;$\par
$(ii)$ on a $p(z)=\log|z|+o(1)$ lorsque $|z|$ tend vers $+\infty;$\par
$(iii)$ il existe $a\in\R$ tel que $p\geqslant a$ dans $K$, et que l'ensemble des $x\in K$ tels que $p(x)>a$ soit polaire.\par
De plus, lorsque ces propriétés sont satisfaites, on a $v(K)=a$.\par}

Il résulte de B.2 et de la prop. B.4.1 que le potentiel $p_K$ de la mesure d'équilibre de $K$ possède ces propriétés, avec $a=v(K)$. Inversement, soit $p:\C\to\overline{\R}$
une fonction  possédant ces propriétés. D'après la prop. B.2.3, elle est de la forme $p_\nu$, où $\nu$ est une mesure de Radon positive dans $\C$, de masse $1$
et à support dans $K$. On déduit de la relation (B.4.3) que l'on a $I(\nu)\geqslant a$, et donc $I(\nu)>-\infty$. Compte tenu des propositions B.3.1 et B.4.1.b), la fonction $p_K$ est $\nu$-presque partout égale à $v(K)$ dans $K$, d'où $\int_K p_K(x)\nu(x)=v(K)$. D'autre part, il résulte de la prop. B.3.1 et de la condition $(iii)$ que la fonction $p_\nu=p$ est $\mu_K$-presque partout égale à $a$ dans $K$, d'où $\int_K p_\nu(y)\mu_K(y)=a$. Or, d'après le théorème de Lebesgue-Fubini ([INT], chap V, \S\ 8, prop.~5), on a
$$\int_Kp_\nu(y)\mu_K(y)=\iint_{K\times K}\log|x-y|\nu(x)\mu_K(y)=\int_\C p_K(x)\nu(x),\leqno({\rm B}.5.2)$$
d'où $v(K)=a$. La prop. B.4.2 implique alors que $\nu$ est la mesure d'équilibre de $K$, et $p$ est son potentiel.

\medskip

\noindent {\bf Corollaire $\bf B.5.3.$} {\it Soit $h:\C\to\R$ une fonction continue, nulle dans $K$ et harmonique dans $\C \sm K$, telle que $h(z)=\log|z|+O(1)$
lorsque $|z|$ tend vers $+\infty$. Alors $\log|z|-h(z)$ a une limite $a$ lorsque $|z|$ tend vers $+\infty$. On a $v(K)=a$ et $p_K=h+a$.} \par

Soit $r>0$ un rayon  tel que $K$ soit contenu dans le disque  $D(0,r)$. La fonction $z\mapsto \log|z^{-1}|-h(z^{-1})$ est alors harmonique et bornée dans le disque ouvert épointé de centre $0$ et rayon $r^{-1}$. Par suite, elle  se prolonge continûment en $0$. D'où l'existence de la limite $a$.

Si $r$ est assez grand, $h$ est positive sur le cercle $C(0,r)$. Le principe du maximum appliqué à $-h$ implique que $h$ est positive dans l'ouvert $D(0,r)\sm K$.
Ainsi $h$ est positive en tout point de $\C$. Sa valeur en un point de $K$ est nulle, donc majorée par la valeur moyenne  de $h$ dans tout disque fermé centré en ce point. Comme par ailleurs $h$ est harmonique dans $\C\sm K$, elle est sous-harmonique dans $\C$ (cf. B.1). Alors $p=h+a$ est sous-harmonique dans $\C$,
égale à $a$ dans $K$ et harmonique dans $\C\sm K$, et l'on a $p(z)=\log|z|+o(1)$ lorsque $|z|$ tend vers $+\infty$. Il résulte de la prop. B.5.1 que l'on a
$v(K)=a$ et $p_K=p$. 

\bigskip

\noindent{\bf B.6. Potentiel de la mesure d'équilibre d'une réunion finie d'intervalles compacts réels}

\smallskip

Soient $g$ un entier naturel et 
$$a_0<b_0<\ldots<a_g<b_g\leqno({\rm B}.6.1)$$
une suite finie strictement croissante de $2g+2$ nombres réels. Comme en $2.3$, posons $E_j=[a_j,b_j]$, notons $E$ la réunion des $E_j$ et posons
$$D(x)=\prod_{j=0}^g(x-a_j)(x-b_j)\leqno({\rm B}.6.2).$$
On a $D(x)\leqslant 0$ pour $x\in E$ et $D(x)>0$ pour $x\in\R\sm E$.\par

D'après $2.3$, il existe un unique polynôme $R\in\R[x]$ unitaire de degré $g$ tel que 
$$\int_{b_{j-1}}^{a_j}{R(x)\over\sqrt{D(x)}}dx=0\qquad\hbox{pour $1\leqslant j\leqslant g$.}\leqno({\rm B}.6.3)$$
Ce polynôme a une unique racine simple dans chaque intervalle $]b_{j-1},a_j[$, où $1\leqslant j\leqslant g$, et aucune autre racine dans $\C$ (cf. 2.3.4).

Il résulte de la démonstration de la prop. 2.3.7 qu'il existe une unique fonction holomorphe $f$ dans $\C\sm E$ telle que $f(x)={R(x)\over\sqrt{D(x)}}$
pour $x\in \ ]b_g,+\infty[$. On a $f^2(z)={R^2(z)\over D(z)}$ et $f(\overline{z})=\overline{f(z)}$ pour tout $z\in\C\sm E$, et $f(z)=z^{-1}+O(z^{-2})$ lorsque 
$|z|$ tend vers $+\infty$.

\medskip

\noindent {\bf Proposition $\bf B.6.4.$} {\it $a)$ Soit $x\in \ ]b_{j-1},a_j[$, où $1\leqslant j\leqslant g$. On a $f(x)=(-1)^{g-j+1}{R(x)\over\sqrt{D(x)}}$.
$b)$ Soit $x\in \ ]a_{j},b_j[$, où $0\leqslant j\leqslant g$. Lorsque $z\in\C$ avec $\Im(z)>0$ tend vers $x$, $f(z)$ tend vers
 $(-1)^{g-j}{R(x)\over i\sqrt{-D(x)}}={|R(x)|\over i\sqrt{|D(x)}|}\cdot$}\par
 
 La dernière égalité résulte de $(2.3.5)$. Le demi-plan supérieur fermé, privé des points $a_0,b_0,\ldots,a_g,b_g$, est simplement connexe. La fonction rationnelle $R^2/D$
y possède donc une racine carrée continue qui prolonge $f$. Cela implique l'existence des limites considérées en~$b)$, et le fait qu'elles dépendent continûment de $x$. Tout revient donc dans $a)$ et $b)$ à déterminer un signe, qui en fait ne dépend que des intervalles considérés et non de $x$.

Le germe de $f$ en $b_g$ peut s'écrire ${u(z)\over \sqrt{z-b_g}}$, où $u$ est le germe d'une fonction holomorphe, réelle sur l'axe réel et non nulle en $b_g$, et $\sqrt{z-b_g}$ est celui de la détermination principale de la racine carrée de $z-b_g$. Comme on a $f(x)={R(x)\over\sqrt{D(x)}}>0$ pour  $x\in \ ]b_g,+\infty[$, on a $u(b_g)\in\R_+^*$. Il s'ensuit que, lorsque $z\in\C$ avec $\Im(z)>0$ tend vers un point $x$ de $]a_g,b_g[$ suffisamment proche de $b_g$, $f(z)$ tend vers un élément de $-i\R_+^*$. L'assertion $b)$ s'en déduit pour $j=g$.

Une étude analogue au voisinage de $a_j$ (resp de $b_{j-1}$), pour $1\leqslant j\leqslant g$, montre que l'assertion $b)$ pour l'indice $j$ implique l'assertion $a)$ pour l'indice $j$ (resp. que l'assertion $a)$ pour l'indice $j$ implique l'assertion $b)$ pour l'indice $j-1$). La proposition en résulte donc par une récurrence descendante.\medskip

{\it Remarque} B.6.5. D'une étude analogue en $a_0$, on déduit  que $f(x)=(-1)^{g+1}{R(x)\over\sqrt{D(x)}}$ pour $x\in \ ]-\infty,a_0[$.\par
{\it Remarque} B.6.6. Soit $x\in \
 ]a_{j},b_j[ \ $, où $0\leqslant j\leqslant g$. Lorsque $z\in\C$ avec $\Im(z)<0$ tend vers $x$, $f(z)$ tend vers
 $(-1)^{g-j}{iR(x)\over \sqrt{-D(x)}}={i|R(x)|\over \sqrt{|D(x)}|}\cdot$. Cela résulte de la prop. B.6.4.b), puisque l'on a $f(\overline{z})=\overline{f(z)}$ pour tout $z\in\C\sm E$.
 
 \medskip
 
 La fonction $f$ possède dans le demi-plan supérieur ouvert une primitive holomorphe. Celle-ci se prolonge en une fonction continue $F_+$ dans  le demi-plan supérieur
 fermé, vu la croissance de $f$ au bord. Normalisons la primitive en exigeant que $F_+(b_g)=0$ et posons alors $h_+=\Re(F_+)$. On a
 $h_+(z)=\log|z|+O(1)$ lorsque $|z|$ tend vers $+\infty$, avec $\Im(z)\geqslant 0$. \medskip
 
\noindent {\bf Proposition $\bf B.6.7.$} {\it La fonction $h_+$ est nulle dans $E$.}\par

Il résulte de la prop. B.6.4.b) que, si $x$ et $y$ sont deux points d'un intervalle  $]a_j,b_j[$, on a
$$h_+(y)-h_+(x)=\Re(\int_x^y{|R(t)|\over i\sqrt{|D(t)|}}dt)=0,$$
donc $h_+$ est constante dans chacun des intervalles $[a_j,b_j]$, où $0\leqslant j\leqslant g$. Par ailleurs, il résulte  de la prop. B.6.4.a) et de la relation (B.6.3) 
que l'on a $h_+(b_{j-1})=h_+(a_j)$ pour $1\leqslant j\leqslant g$. Comme $h_+(b_g)=0$ par hypothèse, $h_+$ est nulle dans chacun des intervalles $[a_j,b_j]$.

\medskip

On définit de manière analogue une fonction continue $F_-$ dans  le demi-plan inférieur fermé, nulle en $b_g$, qui est une primitive de $f$ dans  le demi-plan inférieur 
ouvert. On a en fait $F_-(z)=\overline{F_+(\overline z)}$. On pose
 $h_-=\Re(F_-)$. On a donc $h_-(z)=h_+(\overline{z})$ pour tout $z\in\C$ tel que $\Im(z)\leqslant 0$.\par
 
 \medskip
 
 Les fonctions $h_+$ et $h_-$ coïncident dans $\R$. Notons $h$ la fonction de $\C$ dans $\R$ qui les prolonge. Elle est continue, nulle dans $E$ et harmonique dans 
 $\C\sm E$, puisqu'au voisinage de chaque point de $\C\sm E$ elle peut s'écrire comme la partie réelle d'une primitive de $f$. On a 
 $h(z)=\log|z|+O(1)$ lorsque $|z|$ tend vers $+\infty$. On déduit alors du cor. B.5.2 :
 
 \medskip
 
 \noindent {\bf Proposition $\bf B.6.8.$} {\it On a $h(z)=\log|z|-v(E)+o(1)$ lorsque $|z|$ tend vers $+\infty$, et $h+v(E)$ est le potentiel de la mesure d'équilibre de $E$.}
 
 \medskip

On déduit de la prop. B.6.8 l'expression suivante du logarithme de la capacité de $E$:
$$v(E)=\lim_{x\to+\infty} (\log x-\int_{b_g}^x{R(t)\over\sqrt{D(t)}}dt).\leqno({\rm B}.6.9)$$
 Compte tenu de la remarque B.6.5, on a aussi :
$$v(E)=\lim_{x\to-\infty} (\log |x|-(-1)^{g}\int_x^{a_0}{R(t)\over\sqrt{D(t)}}dt).\leqno({\rm B}.6.10)$$

 \newpage

\noindent{\bf B.7. Mesure d'équilibre d'une réunion finie d'intervalles compacts réels}\smallskip

Conservons les notations de B.6.

\medskip

 \noindent {\bf Théorème $\bf B.7.1.$} {\it La mesure d'équilibre $\mu_E$ de $E$ est la mesure de densité ${|R(x)|\over\pi\sqrt{|D(x)|}}$ par rapport à la mesure de Lebesgue sur $E$.}\par
 
 Nous avons vu dans la prop. B.6.8 que le potentiel $p_E$ de $\mu_E$ est égal à $h+v(E)$, où $h$ est une fonction continue dans $\C$ définie en B.6. 
 Il résulte alors de B.2 que $2\pi\mu_E$ est le laplacien de $p_E$ au sens des distributions; c'est donc aussi le laplacien de $h$. La mesure $\mu_E$ est diffuse (cf. A.4) et son support est contenu dans $E$. Pour démontrer la prop. B.7.1, il nous suffit donc de démontrer que, pour toute fonction  $\varphi:E\to\C$ dont le support est contenu dans l'un des intervalles ouverts $]a_j,b_j[$ et qui est de classe $C^\infty$, on a
 $$2\pi\int_E\varphi(x)\mu_E(x)dx=2\int_{a_j}^{b_j}\varphi(x){|R(x)|\over\sqrt{|D(x)|}}dx.\leqno({\rm B}.7.2)$$
 Choisissons un intervalle $[c,d]\subset \ ]a_j,b_j[$ tel que le support de $\varphi$ soit contenu dans $]c,d[$, puis un nombre réel $r>0$. Il existe une fonction $\psi:\C\to\C$ de classe $C^\infty$ prolongeant $\varphi$ dont le support est contenu dans l'intérieur du rectangle $K= [c,d]  + i [-r,r]$ \ . Comme $2\pi\mu_E=\Delta h$, le premier membre de (B.7.2) est  égal à 
 $$\iint_{\C} \Delta \psi(x+iy)h(x+iy)dxdy=\iint_{K} 2ih{\partial^2\psi\over\partial z\partial\overline{z}}dz\wedge d\overline{z}.\leqno({\rm B}.7.3)$$
 Pour $0<\varepsilon<r$, posons \ $K_\varepsilon^+=  [c,d]  +  i [\varepsilon,r]$ \ et \
  $K_\varepsilon^-=  [c,d] +  i [-r,-\varepsilon]$.
  
  \smallskip
 
 \n Dans le demi-plan supérieur ouvert, on a, avec les notations de B.6, $2{\partial h\over \partial z}={\partial F_+\over \partial z}=f$ et ${\partial^2 h\over \partial z\partial\overline{z}}=0$. Il résulte du théorème de Stokes que l'intégrale $\iint_{K_\varepsilon^+} 2ih{\partial^2\psi\over\partial z\partial\overline{z}}dz\wedge d\overline{z}$ est égale~à 
 $$\int_{\partial K_\varepsilon^+} 2i(\psi {\partial h\over\partial z}dz+h{\partial\psi\over\partial\overline{z}}d\overline{z})=\int_{c+i\varepsilon}^{d+i\varepsilon}i(\psi fdz+2h
 {\partial\psi\over\partial\overline{z}}d\overline{z}).\leqno ({\rm B}.7.4)$$
 Lorsque $\varepsilon$ tend vers $0$, cette intégrale tend vers 
 $$\int_c^d \varphi(x){|R(x)|\over\sqrt{|D(x)|}}dx\leqno ({\rm B}.7.5)$$
 d'après les propositions B.6.4.b) et B.6.7. On démontre de manière analogue que l'intégrale $\iint_{K_\varepsilon^-} 2ih{\partial^2\psi\over\partial z\partial\overline{z}}dz\wedge d\overline{z}$ tend vers cette même limite. La formule (B.7.2) en résulte.

 \newpage


  \section*{ \hspace{55mm} Références}

\bgroup
\leftskip4.2\parindent\parindent0pt
\def\mybibitem[#1]{\removelastskip\leavevmode\hskip-\leftskip{\hbox to\leftskip{[#1]\hss}}\ignorespaces}

\medskip
\mybibitem
[Ab 26] N.H. Abel, {\it Sur l'intégration de la forme différentielle $\rho dx/\sqrt{R}$, $R$ et $\rho$ étant des fonctions entières}\footnote{ Le texte original, écrit en français, est celui reproduit dans les {\it Oeuvres}; celui du journal de Crelle est une traduction en allemand, due probablement à Crelle.},  J. Crelle 1 (1826), 105-144 ($=$ {\it Oe}, XI).

\mybibitem
[ACZ 18] Y. André, P. Corvaja \& U. Zannier, {\it The Betti map associated to a section of an abelian scheme} (with an Appendix by Z. Gao), arXiv:1802.03204[mathAG], à paraître.

\mybibitem
[Ak 30] N.I. Akhiezer (= Achieser = Achyeser), {\it Sur les polynomes de Tchebyscheff pour deux segments}, C.R.A.S. 191 (1930), 754-756.

\mybibitem
[Ak 32]  ------ , {\it Über einige Funktionen welche in zwei gegebenen Intervallen am wenigsten von Null abweichen} I, Izv. Akad. Nauk SSSR (1932), 1163-1202; II, III, {\it ibid.}  (1933), 499-536.

\mybibitem
[A IV] N. Bourbaki, {\it Algèbre. Chapitres} IV-VII, Masson, 1981; traduction anglaise, {\it Algebra} II, Springer-Verlag, 2003.

\mybibitem
[BG 17] A. Bogatyrëv \& O.A. Grigoriev, {\it Closed formula for the capacity of several aligned segments}, Proc. Inst. Steklov Math. 298 (2017), 60-67.

\mybibitem
[Bo 99] A. Bogatyrëv, {\it Effective computations of Chebyshev polynomials for several intervals} (en russe), Math. Sbornik
190 (1999), 15-50; traduction anglaise, Sbornik, Mathematics 190 (1999), 1571-1605.

\mybibitem
[Bo 05]  ------ , {\it Extremal Polynomials and Riemann Surfaces} (en russe), MCCME, Moscou, 2005;
traduction anglaise, Springer-Verlag, 2012.

\mybibitem
[Ca 69] D. Cantor, {\it On approximation by polynomials with algebraic integer coefficients}, Proc. Symposia Pure Math. 12, AMS, 1969, 1-13.

\mybibitem
[Ca 80]  ------ , {\it On an extension of the definition of transfinite diameter and some applications}, J. Crelle 316 (1980), 160-207.

\mybibitem
[Ch 54] P.L. Tchebychef (= Chebyshev), {\it Théorie des mécanismes connus sous le nom de parallélogrammes}, Mém. Acad. Sci. Pétersb. 7 (1854), 539-568 ($=$ {\it Oe} I, 111-143).

\mybibitem
[Ch 58] G. Choquet, {\it Capacitabilité en potentiel logarithmique}, Bull. Acad. Royale Belg. Cl. Sci., 44 (1958), 321-326.

\mybibitem
[Fe 23] M. Fekete, {\it \"{U}ber die Verteilung der Wurzeln bei gewissen algebraischen Gleichungen mit ganzzahligen Koeffizienten}, Math. Zeit. 17 (1923), 228-249.

\mybibitem
[FRV] N. Bourbaki, {\it Fascicule de Résultats des Variétés}, Hermann, Paris, 1971.

\mybibitem
[FS 55] M. Fekete \& G. Szeg\H o, {\it On algebraic equations with integral coefficients whose roots belong to a given point set}, Math. Zeit. 63 (1955), 158-172.

\mybibitem
[Go 03] R. Godement, {\it Analyse Mathématique} IV : {\it  Intégration et théorie spectrale, analyse harmonique, le jardin des délices modulaires}, Springer-Verlag, 2003; traduction anglaise, {\it Analysis} IV, Springer-Verlag, 2015

\mybibitem
[INT] N. Bourbaki, {\it Intégration. Chapitres} I-V, seconde édition, Hermann, Paris, 1967; traduction anglaise, {\it Integration} I, Springer-Verlag, 2004.

\mybibitem
[Kr 57] L. Kronecker, {\it Zwei Sätze über Gleichungen mit ganzzahligen Coefficienten}, J. Crelle 53 (1857), 173-175 (= {\it Werke} I, 103-108).

\mybibitem
[La 16] B. Lawrence, {\it A density result for real hyperelliptic curves}, C.R.A.S. 354 (2016), 1219-1224. 

\mybibitem
[LSN 17] J. Liesen, O. Sète \& M.M.S. Nasser, {\it Fast and accurate computation of the logarithmic capacity of compact sets}, Comp. Methods Funct. Theory 17 (2017), 689-713.

\mybibitem
[Pe 90] F. Peherstorfer, {\it On Bernstein-Szeg\H o orthogonal polynomials on several intervals}, Siam J. Math. Anal. 21 (1990), 461-482.

\mybibitem
[PS 99] F. Peherstorfer \& K. Schiefermayr, {\it Description of extremal polynomials on several intervals and their computation $I, II$}, Acta Math. Hung. 83 (1999), 71-102 \& 103-128.

\mybibitem
[Ra 95] T. Ransford, {\it Potential Theory in the Complex Plane}, London Math. Soc., Students Texts 28, 1995.

\mybibitem
[RR 07] T. Ransford \& J. Rostand, {\it Computation of capacity}, Math. Comp. 76 (2007), 1499-1520.

\mybibitem
[Ro 62] R.M. Robinson, {\it Intervals containing infinitely many sets of conjugate algebraic integers}, in {\it Studies in Mathematical Analysis and Related Topics $:$ Essays in honor of George P\'{o}lya}, Stanford 1962, 305-315.

\mybibitem
[Ro 64]  ------ , {\it Conjugate algebraic integers in real point sets}, Math. Zeit. 84 (1964), 415-427.

\mybibitem
[Ro 69]  ------ , {\it Conjugate algebraic integers on a circle}, Math. Zeit. 110 (1969), 41-51.

\mybibitem
[Ru 14] R. Rumely, {\it Capacity Theory with Local Rationality - The strong Fekete-Szeg\H o theorem on curves},
A.M.S. Surveys 193, 2013.

\mybibitem
[Sc 18] I. Schur, {\it Über die Verteilung der Wurzeln bei gewissen algebraischen Gleichungen mit ganzzahligen Koeffizienten}, Math. Zeit. 1 (1918), 377-402 (= {\it Gesam. Abh.} II, 32). 

\mybibitem
[Se 97] J-P. Serre, {\it Répartition asymptotique des valeurs propres de l'opérateur de Hecke $T_p$}, J.A.M.S. 10 (1997), 75-102 ($=$ {\it Oe.} IV, 170).

\mybibitem
[Sm 84] C. Smyth, {\it Totally positive algebraic integers of small trace}, Ann. Inst. Fourier 33 (1984), 1-28.

\mybibitem
[SY 92] M.L. Sodin \& P.M. Yuditskii, {\it Functions which deviate least from zero on closed subsets of the real axis} (en russe), Algebra i Analiz 4 (1992), 1-61; traduction anglaise, St. Petersburg Math. J. 4 (1993), 201-249.

\mybibitem
[St 85] T.J. Stieltjes, {\it Sur quelques théorèmes d'algèbre}, C.R.A.S. 100 (1885), 439-440 ($=$ {\it Oe.} I, 440-441).

\mybibitem
[Sz 24] G. Szeg\H o, {\it Bermerkungen über einer Arbeit von Herrn M. Fekete~$:$\"{U}ber die Verteilung der Wurzeln bei gewissen algebraischen Gleichungen mit ganzzahligen Koeffizienten}, Math. Zeit. 21 (1924), 203-208.
 
\mybibitem
[Ta 69]  J. Tate, {\it Classes d'isogénie des variétés abéliennes sur un corps fini}, Sém. Bourbaki 1968/1969, n° 352  
(= LNM 175 (1971), 95-110 = {\it Coll. Papers} I, 32).

\mybibitem
[Ts 18] M. A. Tsfasman, {\it Serre's theorem and measures corresponding to abelian varieties over finite fields}, à
paraître.
\mybibitem
[Ts 59] M. Tsuji, {\it Potential Theory in Modern Function Theory}, Maruzen, Tokyo, 1959; seconde édition, Chelsea, New York, 1975.

\mybibitem
[TV 97] M.A. Tsfasman \& S.G. Vl\u{a}du\c{t}, {\it Asymptotic properties of zeta functions}, J.~Math. Sci. (New York) 84 (1997), 1445-1467.

\egroup

%
%
%
%
%
%
%
%

\end{document}